\documentclass[12pt,a4paper]{article}

\usepackage{amsmath}
\usepackage{amssymb}
\usepackage{graphics}

\newcommand{\mr}[1]{\mathrm{#1}}
\newcommand{\mi}[1]{\mathit{#1}}

\newcommand{\fmi}[1]{\underline{\mathit{#1}}}
\newcommand{\mg}[1]{{\text{\boldmath$#1$}}}

\newcommand{\pp}{{\scriptscriptstyle +}}
\newcommand{\mm}{{\scriptscriptstyle -}}
\newcommand{\ppmm}{{\scriptscriptstyle \pm}}
\newcommand{\mmpp}{{\scriptscriptstyle \mp}}
\newcommand{\p}{\!+\!}
\newcommand{\m}{\!- \!}
\newcommand{\bu}{\bullet}

\newcommand{\gop}{\mathfrak{p}}

\newcommand{\caA}{\mathcal{A}}
\newcommand{\caB}{\mathcal{B}}

\newcommand{\caL}{\mathcal{L}}

\newcommand{\caO}{\mathcal{O}}

\newcommand{\caR}{\mathcal{R}}

\newcommand{\caU}{\mathcal{U}}

\newcommand{\doA}{\mathbb{A}}

\newcommand{\doC}{\mathbb{C}}

\newcommand{\doN}{\mathbb{N}}

\newcommand{\doQ}{\mathbb{Q}}

\newcommand{\doS}{\mathbb{S}}

\newcommand{\doZ}{\mathbb{Z}}

\newcommand{\doDelta}{\Delta \hspace{-1.6 ex} \Delta}

\begin{document}

\title{Invariants of identity-tangent diffeomorphisms: explicit formulae and effective computation.}
\author{O. Bouillot and J. Ecalle; Orsay, France.}
\maketitle



%

\noindent
{\bf Abstract:} {\it In this short Survey we revisit the subject of local, identity-tangent diffeomorphisms of $\doC$ and their analytic invariants, under two viewpoints: that of explicit expansions, which necessarily involve multitangents and multizetas; and that  of effective computation. Along the way, we stress the difference between the \textup{collectors} (pre-invariant but of one piece) and the \textup{connectors} (invariant but mutually unrelated). 
We also attempt to streamline the nomenclature and notations. 
}

\tableofcontents

%


\section{Setting and notations.}


\subsection{Classical results.}
We shall be concerned here with {\it local\,\footnote{
i.e. {\it analytic germs of \---}} identity-tangent diffeomorphisms} of $\doC$,  or {\it diffeos} for short, with the fixed-point located at $\infty$ for technical convenience: 
\begin{equation}  \label{wa1}
f\;\; : \quad z \mapsto z+\sum_{1\leq s} f_s\; z^{1-s}     \hspace{7.ex}  a_s\in \doC
\end{equation}
Unless $f$ be the identity map, we can always subject it to an analytic (resp. formal) conjugation 
$f \mapsto  f_1=h\circ f \circ h^{-1}$, followed if necessary by an elementary ramification $\big(f_1(z^{1/p})\big)^p$, so as to give it the following {\it prepared} resp. {\it normal} form:
\begin{eqnarray} \label{prep}
f_{\mr{prep}}: \;\; z& \mapsto & z+1-\rho z^{-1}+\sum_{2< s_0\leq s} f_{[s]}\; z^{1-s}  \hspace{5.ex} (s\in \frac{1}{p}\doN^\ast)
\\  \label{norm}
f_{\mr{norm}}:\;\; z& \mapsto & z+1-\rho z^{-1}
 \end{eqnarray}
where $s_0$ may be chosen as large as one wishes.

The {\it tangency order} $p$ and {\it iterative residue} $\rho$ are the only {\it formal invariants} of identity-tangent diffeos. But our diffeos also possess countably many (independent) scalar {\it analytic invariants}, which are best defined as the Fourier coefficients of the so-called {\it connectors}.\footnote{In the context of identity-tangent diffeos, the connectors are sometimes referred to as {\it horn maps}, but the notion is more general: in resurgent analysis (see \S1.2 {\it infra}) the connectors are the operators that take us from on sectorial model to the next.
}
 These 
are pairs of germs of 1-periodic analytic mappings $\mg{\pi}=(\mg{\pi}_{\mr{no}},\mg{\pi}_{\mr{so}}) $ on upper/lower half-planes $\pm \Im(z)\gg1$. There are $p$ such pairs, corresponding to the $p$-fold ramification of $z$ in (\ref{prep}).
Here, {\it no} and {\it so} stand for {\it north} and {\it south}, i.e. the upper and lower half-planes.

We shall throughout prioritise the {\it standard case} $p=1$ , $\rho=0$, i.e. focus on diffeos of the form:
\begin{equation}  \label{wa2}
f:=l \circ g \quad \mi{with} \quad
 l\,:\,= z\mapsto z+1 \quad \mi{and} \quad g\,:\, z\mapsto z+\sum_{3\leq s}g_s\, z^{1-s}
\end{equation}
and merely sketch the (minor) changes required to cover the general case.

Any standard $f$ possesses two well-defined, mutually inverse so-called {\it iterators}, to wit $f^\ast_{\pm}$ (direct iterator)
and  $^{\ast\!\!}f_{\pm}$ (reciprocal iterator), defined on U-shaped domains\footnote{ $f^\ast_{+}$ and $^{\ast\!\!}f_{+}$ 
are defined on a west-north-south domain, while $f^\ast_{-}$ and $^{\ast\!\!}f_{-}$ are defined on an east-north-south domain.
}
 by the limits:
\begin{equation}    \label{ekam1}
f^\ast_{\ppmm}(z) =\lim_{k\rightarrow \ppmm \infty} l^{-k}\circ f^{k}
\quad;\quad
{^\ast\!}{f_{\ppmm}}(z) =\lim_{k\rightarrow \ppmm \infty} f^{-k}\circ l^{k}
\end{equation}
The connectors $\mg{\pi}^{\pm1}$, with their northern and southern components, are then defined on 
$\pm\Im(z)\gg 1$ by:
\begin{equation}   \label{ekam2}
\mg{\pi}:=f^\ast_\pp\circ {^\ast\!}{f_{\mm}}
\quad;\quad
\mg{\pi}^{-1}:=f^\ast_\mm\circ {^\ast\!}{f_{\pp}}
\end{equation}
For reasons that will soon become apparent, we must also consider the infinitesimal
generators $f_\ast$ and $\mg{\pi}_\ast$ of $f$ and $\mg{\pi}$. These are formal, generically divergent power resp. Fourier series. Of course, $\mg{\pi}_\ast$ is not constructed directly from  $\mg{\pi}$, but via its northern and southern components.   
We thus have the three pairs:
\begin{equation}   \label{ekam3}
\mg{\pi}:=(\mg{\pi}_{\mr{no}},\mg{\pi}_{\mr{so}})
\quad;\quad
\mg{\pi}^{-1}:=(\mg{\pi}^{-1}_{\mr{no}},\mg{\pi}^{-1}_{\mr{so}})
\quad;\quad
\mg{\pi}_{\ast}:=(\mg{\pi}_{\ast\mr{no}},\mg{\pi}_{\ast\mr{so}})
\end{equation}
along with the relations
\begin{eqnarray}   \label{ekam4}
 f(z) &=& \exp\!\big( f_\ast(z)\,\partial_z \big).\, z  
 \hspace{10. ex}  \big(\,f_\ast \,\partial_z f^\ast \equiv1\,\big)
\\ [1.ex]   \label{ekam5}
\mg{\pi}^{\pm1}_{\mr{no}}(z)&=& \exp\!\big(\,\pm\,\mg{\pi}_{\ast\mr{no}}(z)\, \partial_z  \big).\,z
\\ [1.ex]   \label{ekam6}
\mg{\pi}^{\pm1}_{\mr{so}}(z)&=& \exp\!\big(\,\pm\,\mg{\pi}_{\ast\mr{so}}(z)\, \partial_z  \big).\,z
\end{eqnarray}
In   (\ref{ekam4})  $f^\ast$ and $^{\ast\!\!}f$ denote of course the {\it formal iterators}, i.e. the power series solutions of the equations
\begin{eqnarray}    \label{iterat1}
f^\ast\circ f =\,l\circ f^\ast  \quad &\mi{with}& \quad   f^\ast(z)=z+o(1)
\\    \label{iterat2}
f\,\circ\, ^{\ast\!\!}f =\, ^{\ast\!\!}f \circ l  \quad &\mi{with}& \quad  ^{\ast\!\!}f(z)=z+o(1)
\end{eqnarray}
normalised by the condition of carrying no constant term. Anticipating on the sequel, here is how the scalar invariants can be read off the Fourier expansions of the connectors:
\begin{eqnarray}   \label{ekam7}
\mg{\pi}_{\mr{no}}(z)=z+\sum_{\omega\in\Omega^-} A^{+}_{\omega}\,e^{-\omega\,z}
&;&
\mg{\pi}_{\mr{so}}(z)=z+\sum_{\omega\in\Omega^+} A^{-}_{\omega}\,e^{-\omega\,z}
\\   \label{ekam8}
\mg{\pi}^{-1}_{\mr{no}}(z)\!=z+\,\sum_{\omega\in\Omega^-} A^{-}_{\omega}\,e^{-\omega\,z}
&;&
\mg{\pi}^{-1}_{\mr{so}}(z)=z+\sum_{\omega\in\Omega^+} A^{+}_{\omega}\,e^{-\omega\,z}
\\   \label{ekam9}
\mg{\pi}_{\ast\mr{no}}(z)=\!+2\pi i \!\!\sum_{\omega\in\Omega^-} A_{\omega}\,e^{-\omega\,z}
&;&
\mg{\pi}_{\ast\mr{so}}(z)=\!-2\pi i \!\! \sum_{\omega\in\Omega^+} A_{\omega}\,e^{-\omega\,z}
\end{eqnarray}
Pay attention to the altered position of $\pm$ in $\ref{ekam7}$ and  $\ref{ekam8}$; the reasons for this apparent incoherence shall become clear in due course. The indices $\omega$ run through $\Omega:=2\pi i\doZ^\ast$ or
 $\Omega^\pm:=\,\pm 2\pi i\, \doN^\ast$, and each of the three systems
 \begin{equation}  \label{wa3}
 \{ A^{+}_\omega\,,\, \omega\in \Omega\}\quad\quad, \quad\quad
 \{ A^{-}_\omega\,,\, \omega\in \Omega\}\quad\quad , \quad\quad
  \{ A_\omega\,,\, \omega\in \Omega\}\quad \quad
 \end{equation}
constitutes a {\it free} and {\it complete} system of analytic invariants.\footnote{ With the minor qualifier that, under a conjugation by a shift $h$ of the form $l^{\alpha}(z):=z+\alpha$, the periodic germs 
$\mg{\pi}^{\pm}$
also undergo conjugation by the same shift, with obvious repercussions for their Fourier coefficients.
}
%
\subsection{Reminder about resurgent functions.}
We will have to be content here with a very sketchy presentation. The algebra of {\it resurgent fonctions} admits three different realisations or models:
\\
(i) the {\it formal model}, consisting of formal power series $\tilde{\varphi}(z)$ of $z^{-1}$ or of more general {\it transseries};
\\
(ii) the {\it convolutive model}, consisting of microfunctions\footnote{
i.e. minor-major pairs 
  $(\hat{\varphi}(\zeta),\check{\varphi}(\zeta))$. The {\it majors} are defined up to regular germs at the origin, 
and the {\it minors} are related to them under 
$ 2\pi i\,\hat{\varphi}(\zeta)\equiv \check{\varphi}(\zeta\,e^{-\pi i})- \check{\varphi}(\zeta\,e^{+\pi i}) $
for $\zeta\sim0$. In the present paper, we shall almost entirely dispense with majors, since we shall mostly 
be dealing with so-called {\it integrable} microfunctions, whose minors carry the whole information.
}
 at $\zeta=0$, whose {\it majors}  $\check{\varphi}(\zeta)$ 
are constraint-free at the origin but whose
{\it minors} $\hat{\varphi}(\zeta)$ have the property of endless continuation\footnote{laterally along any {\it finite and finitely punctured} broken lines.} and exponential growth;\footnote{i.e. {\it at most exponential}, along {\it infinite but finitely punctured broken lines}, with a suitable uniformity condition. }
\\
(iii) the {\it geometric model(s)}, consisting of analytic germs $\varphi_\theta(z)$ defined on sectorial neighbourhoods of 
$\infty$ of bisectrix $\mi{arg}( z^{-1})=\theta$ and aperture at least $\pi$. 

The natural algebra product in the $z$-models (i) and (iii) is of course multiplication. In the $\zeta$-model (ii) it is convolution, defined first {\it locally}\footnote{When  the {\it minors} $\hat{\varphi}$  are not integrable at the origin, one must modify the definition and
draw in the {\it majors} $\check{\varphi}$.  }  by
\begin{equation}   \label{wa4}
(\hat{\varphi}_1 \ast \hat{\varphi}_2)(\zeta):=
\int_0^\zeta \hat{\varphi}_1(\zeta_1)\,\hat{\varphi}_2(\zeta-\zeta_1)\,d\zeta_1
\hspace{6.ex} (\zeta\sim 0)
\end{equation}
and then {\it in the large} by analytic continuation.

In practice, one starts with elements $\tilde{\varphi}$ of model (i) obtained as formal solutions of differential or functional equations, and the aim is to resum them, i.e. to go to model (iii). Generally speaking, this is possible only over the detour through model (ii), with the {\it formal Borel tranform} $\caB$
\begin{equation}  \label{wa5}
z^{-\sigma}\mapsto \frac{\zeta^{\sigma-1}}{\Gamma(\sigma)}
\quad ; \quad 
(\partial_\sigma)^n z^{-\sigma}\mapsto 
(\partial_\sigma)^n \frac{\zeta^{\sigma-1}}{\Gamma(\sigma)}
\quad ; \quad \mi{etc}
\end{equation}
taking us from (i) to (ii), and the {\it polarised Laplace transform} $\caL_\theta$
\begin{equation}  \label{wa6}
\varphi_\theta(z)=\int_{\arg(\zeta)=\theta}\hat{\varphi}(\zeta)\,e^{-\zeta z}\, d\zeta
\end{equation}
taking us from (ii) to (iii).

The most outstanding feature of the resurgence algebras is the existence on them of a rich array of so-called {\it alien operators} $\Delta_\omega$ and  $\Delta^\pm_\omega$, with indices $\omega$ running through 
$\doC_\bu:=\widetilde{\doC - \{0\}}$. These operators act on all three models\footnote{with the same symbols doing service in all three, since no confusion is possible.}, but are first defined in the convolutive model, where they have the effect of measuring the singularities of the (often highly ramified) minors $\hat{\varphi}$ at or rather {\it over} $\omega$.
Here is how they act:
\begin{eqnarray}   \label{ekam11}
(\Delta_\omega\hat{\varphi})(\zeta)&:=& \sum_{\epsilon_1,...,\epsilon_r}
\frac{\epsilon_r}{2\pi i}\, \lambda_{\epsilon_1,...,\epsilon_{r-1}}\; 
\hat{\varphi}^{({\epsilon_1 \atop \omega_1 }{,..., \atop ,..., }{\epsilon_1 \atop \epsilon_r })}(\omega+\zeta)
\\   \label{ekam12}
(\Delta^\pm_\omega\hat{\varphi})(\zeta)&:=& \sum_{\epsilon_1,...,\epsilon_r}
\pm\,\epsilon_r\, \lambda^{\pm}_{\epsilon_1,...,\epsilon_{r-1}}\; 
\hat{\varphi}^{({\epsilon_1 \atop \omega_1 }{,..., \atop ,..., }{\epsilon_1 \atop \epsilon_r })}(\omega+\zeta)
\end{eqnarray}
with signs $\epsilon_j\in \{+,-\}$,
with weights $\lambda_\bu$, $\lambda^{+}_\bu$, $\lambda^{-}_\bu$ defined by
\begin{eqnarray}   \label{ekam13}
\lambda_{\epsilon_1,...,\epsilon_{r-1}} &:=& \frac{p!\,q!}{r!}\quad \mi{with}
\quad p:=\sum_{\epsilon_i=+}\,1\;,\; q:=\sum_{\epsilon_i=-}\,1\;\; 
\\   \label{ekam14}
\lambda^\epsilon_{\epsilon_1,...,\epsilon_{r-1}} &:=& 1 \quad \mi{if} \quad \epsilon_1=\dots=\epsilon_{r-1}=\epsilon
\\   \nonumber                                                     &:=& 0 \quad \mi{otherwise}
\end{eqnarray}
and with 
$\hat{\varphi}^{({\epsilon_1 \atop \omega_1 }{,..., \atop ,..., }{\epsilon_1 \atop \epsilon_r })}(\omega+\zeta)$
denoting the analytic continuation of $\hat{\varphi}$ from $\zeta$ to $\omega+\zeta$ under right (resp. left) circumvention of each intervening singularity $\omega_j$ if $\epsilon_j=+$ (resp. $\epsilon_j=-$).  We start of course with a point $\zeta$
close enough to 0 on the axis $\mi{arg}(\zeta)=\mi{arg}(\omega)$, and extend the definition in the large by analytic continuation. The operators $\Delta_\omega$ are {\it derivations}. Thus, in the convolutive and formal models, the identities hold:
\begin{eqnarray}   \label{ekam10}
\Delta_\omega(\hat{\varphi}_1\ast \hat{\varphi}_2) &= &
  \Delta_\omega(\hat{\varphi}_1)\ast \hat{\varphi}_2 
+\hat{\varphi}_1\ast \Delta_\omega(\hat{\varphi}_2) 
\\   \label{ekam10bis}
\Delta_\omega(\tilde{\varphi}_1\;.\; \tilde{\varphi}_2) &=& 
  \Delta_\omega(\tilde{\varphi}_1)\;.\; \tilde{\varphi}_2 
+\tilde{\varphi}_1\;.\; \Delta_\omega(\tilde{\varphi}_2) 
\end{eqnarray}
When working in any one of the multiplicative models (formal or geometric), it is often convenient to phase-shift the alien operators, and to set:
\begin{eqnarray}   \label{ekam15}
\doDelta_\omega &:=& e^{-\omega z} \Delta_\omega 
\quad\quad\quad (\;[\partial_z, \doDelta_\omega]\equiv 0\;)
\\ [1. ex]   \label{ekam16}
\doDelta_\omega^{\pm} &:=& e^{-\omega z} \Delta_\omega^{\pm}
\quad\quad\quad (\;[\partial_z, \doDelta^{\pm}_\omega]\equiv 0\;)
\end{eqnarray}
The gain here is that the new operators commute with $\partial_z$. These phase-shifted operators are also the natural ingredients of the {\it axial operators} $\caR_{_\theta}$ and  $\caR_{_\theta}^{\pm}$\;:
\begin{eqnarray}    \label{ekam17ante}
\caR_{_{\theta}} &=&\sum_{\arg(\omega)=\theta}\!\! \doDelta_\omega
\\  \label{ekam17}
\caR_{_\theta}^\pm &=&
 1+\!\!\!\sum_{\arg(\omega)=\theta} \, \doDelta^{\pm}_\omega
 \;=\; 
\exp\big(\,\pm 2\pi i \, \caR_{_{\theta}}\big) 
\end{eqnarray}
which are the key to the {\it axis-crossing} identities\,:
\begin{eqnarray}\label{cross+}  
\varphi_{_{\theta-\epsilon}} = (\caR^{+}_{_\theta}\; \varphi)_{_{\theta+\epsilon}} 
& \quad;\quad &
(\,\varPhi\, . \, \caR^{+}_{_\theta}\,)_{_{\theta-\epsilon}} =
(\,\caR^{+}_{_\theta} . \,\varPhi \, )_{_{\theta+\epsilon}}
\\ [1.ex]\label{cross-}   
\varphi_{_{\theta+\epsilon}} = (\caR^{-}_{_\theta}\, \varphi)_{_{\theta-\epsilon}} 
& \quad;\quad &
(\,\varPhi\, . \, \caR^{-}_{_\theta}\,)_{_{\theta-\epsilon}} =
(\,\caR^{-}_{_\theta} . \,\varPhi \, )_{_{\theta+\epsilon}}
\end{eqnarray}
that connect two sectorial germs $\varphi_{\theta-\epsilon}$ and  $\varphi_{\theta+\epsilon}$ relative to Laplace integration
right and left of any given singularity-carrying axis $\theta$ in the $\zeta$-plane.\footnote{
In (\ref{ekam10}), (\ref{ekam10bis}), $\varphi$ denotes any {\it resurgent function} and $\varPhi$ any {\it
resurgent operator} (such as multiplication or postcomposition by a resurgent function etc).
}
%
\subsection{The Bridge equation.}
The iterator $f^\ast$ and $^{\ast\!\!}f$ characterised by the relations (\ref{iterat1}) and  (\ref{iterat2}) 
verify the following resurgence equations
\begin{eqnarray}   \label{ekam22}
\Delta_\omega\, ^{\ast\!\!}f(z) &=& +A_\omega\, \partial_z\, ^{\ast\!\!}f(z)
\hspace{14.ex} (\forall \omega\in \Omega)
\\ [0.6 ex]   \label{ekam23}
\Delta_\omega\, f^\ast(z) &=& -A_\omega \, e^{-\omega\,(f^\ast(z) -z) }
\hspace{11.ex} (\forall \omega\in \Omega)
\end{eqnarray}
with the very same scalar coefficients $A_\omega$ as in (\ref{ekam9}). For  values of $\omega$
not in $\Omega$, the alien derivatives would be $\equiv 0$. If we now introduce the differential operators:
\begin{equation}    \label{wa10}
\doA_\omega := A_\omega\,e^{-\omega z}\,\partial_z \hspace{11.ex} (\forall \omega\in \Omega)
\end{equation}
the resurgence equations assume the form of the Bridge equation:\footnote{
so-called because it relates {\it ordinary}  and {\it alien} derivatives of one and the same resurgent function. The Bridge equation has in fact much wider applications, and extends, in one form or another, to practically all {\it resonant} local objects, of which {\it identity-tangent diffeos} are but a special case. An entire book [E3] has been devoted to the subject.} 
\begin{eqnarray}     \label{ekam24}
\doDelta_\omega \, ^{\ast\!\!}f(z) &=& +\doA_\omega\, \partial_z\, ^{\ast\!\!}f(z)
\\  [0.9 ex]  \label{ekam24bis}
\doDelta_\omega \, f^{\ast}(z) &=& -(\doA_\omega\,z) \circ f^\ast(z)
\end{eqnarray}
When expressed in terms of the subsitution operators $F^\ast$ and $^{\ast\!\!}F$
associated with ${^\ast\!\!}f, f^\ast$, 
the Bridge equation takes an even more pleasant form
\begin{eqnarray}   \label{ekam26}
\big[  \doDelta_\omega, F^\ast \big] &=& -\,F^\ast\;\doA_\omega
\hspace{10.ex} ( F^\ast\, \varphi := \varphi\circ f^\ast)
\\         \label{ekam27}
\big[ \doDelta_\omega, \,^{\ast\!\!}F \,\big] &=& +\,\doA_\omega\; ^{\ast\!\!}F
\hspace{10.ex} \;(\, ^{\ast\!\!}F\; \varphi := \varphi\circ ^{\ast\!\!}f)
\end{eqnarray}
But whichever variant we may care to consider, the commutation identities $ \big[\doDelta_{\omega_1}, \doA_{\omega_2} \big]=0 $ 
make it easy to iterate the above resurgence equations. Thus from (\ref{ekam24}) we straightaway derive
\begin{equation}   \label{ekam25}
\doDelta_{\omega_r} \dots \doDelta_{\omega_1} \, ^{\ast\!\!}f(z) =
 \doA_{\omega_1}\dots\doA_{\omega_r} \, \partial_z\, ^{\ast\!\!}f(z)\hspace{6.ex}\mi{(order \;reversion!)}
\end{equation}
As a consequence, the effect on $^{\ast\!\!}f$ and $f^\ast$  of the alien operators $\doDelta^{\pm}_\omega$ and of the axial operators $\caR_\theta$ is easy to calculate. 
It is best written in terms of the substitution operators 
${^\ast\!}F$ and $F^\ast$  
associated with ${^\ast\!\!}f, f^\ast$, and results in the so-called
{\it axial} Bridge equation:
\begin{eqnarray}     \label{ekam30pre}   
\caA_\theta&=& \caR_{\theta}   \; - \;  {^{\ast\!\!}F}\; \; \caR_{\theta} \;\; F^\ast
\\ [1.ex]
\label{ekam30}
\caA^{+}_\theta&=& \caR^{+}_\theta\;\; ^{\ast\!\!}F\; \; \caR_\theta^{-}\;\; F^\ast \, =\, 
^{\ast\!\!}F\;\;  \caR^{-}_\theta\;\; F^\ast \;\; \caR_\theta^{+}
\\ [1.ex]
\label{ekam30post}
\caA^{-}_\theta&=& \caR^{-}_\theta\;\; ^{\ast\!\!}F\; \; \caR_\theta^{+}\;\; F^\ast \, =\, 
^{\ast\!\!}F\;\;  \caR^{+}_\theta\;\; F^\ast \;\; \caR_\theta^{-}
\end{eqnarray} 
The axial Bridge equation\footnote{
The {\it singular} is appropriate here since  (\ref{ekam30}) and  (\ref{ekam30post}) are merely exponential variants of
 (\ref{ekam30pre}). The commutation of the three automorphisms
   $\caA^{\pm}_\theta $, $\caR^{\pm}_\theta $, $ ^{\ast\!\!}F\; \; \caR_\theta^{\mp}\;\; F^\ast $  
is itself a consequence of the commutation of the three derivations
  $ \caA_\theta$, $\caR_{\theta} $, ${^{\ast\!\!}F}\; \; \caR_{\theta} \;\; F^\ast $.
}
 involves 
differential (resp. substitution) operators $\caA_\theta$ (resp.  $\caA^\pm_\theta$ ):
\begin{eqnarray}          \label{ekam28}
\caA_\theta &=&  \sum_{\arg(\omega)=\theta} \, \doA_\omega
\\        \label{ekam29}
\caA_{_\theta}^\pm &=&
 1+\sum_{\arg(\omega)=\theta} \, \doA^{\pm}_\omega
 \;=\; 
\exp\big(\,\pm 2\pi i\,\caA_\theta \big) 
\end{eqnarray}
which are simply related to the differential (resp. substitution) operators
$\mg{\Pi}_\ast$ (resp. $\mg{\Pi}^{\pm}$ 
 associated with the connectors of \S1.1:
\begin{eqnarray}   \label{cro1}
\mg{\Pi}_{\mg{no}}\;\;:=\hspace{6.ex} \caA^{+}_{-\frac{\pi}{2}} &\quad;\quad&
\mg{\Pi}_{\mg{so}}\;\;:=\hspace{6.ex}  \caA^{-}_{+\frac{\pi}{2}} 
\\  \label{cro2}
\mg{\Pi}^{-1}_{\mg{no}}\;\;:=\hspace{6.ex} \caA^{-}_{-\frac{\pi}{2}} &\quad;\quad&
\mg{\Pi}^{-1}_{\mg{so}}\;\;:=\hspace{6.ex} \caA^{+}_{+\frac{\pi}{2}} 
\\  \label{cro3}
\mg{\Pi}_{\ast \mg{no}}\; :=\,+2\pi i\,\caA_{-\frac{\pi}{2}} &\quad;\quad&
\mg{\Pi}_{\ast \mg{so}}\; :=\,-2\pi i\,\caA_{+\frac{\pi}{2}} 
\end{eqnarray}
The first identity  (\ref{cro1}) results from applying the direct axis-crossing formula (\ref{cross+})  with $\theta=-\frac{\pi}{2}$
and $\varphi={^{\ast\!\!}f} $ or $\varPhi={^{\ast\!\!}F} $, since $^{\ast\!\!}f_{_{\theta\pm \epsilon}}={^{\ast\!\!}f_\pm}$.
The second identity  (\ref{cro1}) results from applying the inverse axis-crossing formula (\ref{cross-})  with $\theta=+\frac{\pi}{2}$
and $\varphi={^{\ast\!\!}f }$ or $\varPhi={^{\ast\!\!}F} $, since in that case  $^{\ast\!\!}f_{_{\theta\pm\epsilon}}={^{\ast\!\!}f_\mp}$ (inversion!). The identities (\ref{cro2}) and  (\ref{cro2}) immediately follow.
%
 \subsection{Invariants, connectors, collectors.}
Let us survey in one table the main objects introduced so far or yet to come.
\[\begin{array}{ccccccccccccc}
\mi{diffeo} &&  \mi{collectors}  && \mi{connectors} && \mi{invariants}
\\ [1.5 ex]
f=l\circ g &\stackrel{1^{'}}{\longrightarrow} & \mg{\gop}\stackrel{1^{''}}{\longrightarrow} \mg{p} &
\stackrel{1^{'''}}{\longrightarrow} &\mg{\pi}=(\mg{\pi}_{\mr{no}},\mg{\pi}_{\mr{so}})
&\stackrel{1^{''''}}{\longrightarrow} &  \{ A_\omega^{\pm}\}
\\ [1.5 ex]
\downarrow \scriptstyle{2}  & & \uparrow \scriptstyle{4}\hspace{5.ex} &&  
\hspace{5.4 ex}\uparrow \scriptstyle{5_{\mr{no}}} \hspace{2. ex}\uparrow \scriptstyle{5_{\mr{so}}}
& & \uparrow \scriptstyle{6}
\\ [1.5 ex]
g_\ast &\stackrel{3^{'}}{\longrightarrow} &\; \,\mg{\gop}_\ast\! \stackrel{3^{''}}{\longrightarrow} \mg{p}_\ast &
\stackrel{3^{'''}}{\longrightarrow} &\,\mg{\pi}_{\ast}=(\!\mg{\pi}_{\mr{\ast no}},\mg{\pi}_{\mr{\ast so}}\!)
&\stackrel{3^{''''}}{\longrightarrow} & \, \{ A_\omega\}
\end{array}\]

The upper  row carries the objects of direct interest to us,
while the lower row carries their infinitesimal counterparts, which are more in the nature of auxiliary tools.

The first, third and fourth columns carry objects already familiar to us. The second column, however, carries novel, highly interesting objects,  the {\it collectors}, which are very close in a sense to the {\it connectors}, yet should be, for the sake of conceptual cleanness, clearly held apart. The {\it collectors} may assume three distinct forms:
\\
(i) formal series of multitangents, noted $ \mg{\gop}$;
\\
(ii) formal series of monotangents, also noted $ \mg{\gop}$;
\\
(iii) formal power series of $z^{-1}$, noted $ \mg{p}$, and their Borel transforms  $ \mg{\hat{p}}$.

One goes from (i) to (ii) by multitangent reduction as in \S2.3\,; and from (ii) to (iii) by the change 
$\mi{Te}^{s_1}\mapsto z^{-s_1}$.

In any of their three incarnations, the collectors are but a step removed from the invariants. Yet they are not invariant themselves: they depend on the chart in which the diffeo $f$ is taken. In fact, the diffeo can be re-constructed in totality from the collector -- easily so if the collector is taken in the form (i), less easily if it is taken in the more {\it condensed} forms (ii) or (iii).

One last remark is in order here: although we are basically interested in the objects of the upper column, and more specifically in getting from $f$ to its invariants $ \{ A^\pm_\omega\}$, we shall see that the most advantageous route is not the straight path through the arrows
 $1,1^{\prime},1^{\prime\prime},1^{\prime\prime\prime\prime}$, but any of the roundabout paths that start with $2$ and $3^\prime$: these indirect routes are much more economical in terms of calculations and also more respectful of the underlying symmetries and parities.



%
 \subsection{The reverse problem: canonical synthesis.}
It can be shown that any pair $\mg{\pi}=(\mg{\pi}_{\mr{no}}, \mg{\pi}_{\mr{so}})$ is the connector of some standard diffeo $f$. This raises the problem of {\it synthesis}: how to reconstitute some $f$ from a given set of invariants? And how to produce a canonical $f$ among all possible choices? A semi-canonical synthesis was sketched in [E2] and a fully canonical one was constructed in [E4]. The latter depends 
on a single parameter $c$ whose real part must be chosen large enough.\footnote{ 
 Synthesis cannot be {\it absolute}, i.e. parameter-free. } 
 The construction produces a canonical
$f_c:=^{\ast\!\!}f_c\circ l\circ f_c^\ast $ from its iterator $f_c^\ast$, which in turn is explictely
 given, in operator form, by the formula
\begin{equation}        \label{ekam34}
F_c^\ast:=1 +\sum_r\sum_{\omega_i \in \Omega} (-1)^r\; {\caU\!{e}}_c^{\omega_1,\omega_2,...,\omega_r}(z)\;
\doA_{\omega_r} \dots \doA_{\omega_2}\,\doA_{\omega_1}
\end{equation}
with a careful re-arrangement of the terms\footnote{ known as arborification-coarborification.} necessary to ensure convergence. 
The two ingredients in  (\ref{ekam34}) are the invariants $\doA_\omega$ taken in operator form (\ref{wa10}), and some special resurgence monomials $ {\caU\!{e}}_c^{\mg{\omega}}(z)$ defined by
\begin{equation}        \label{ekam35}
{\caU\!{e}}_c^{\mg{\omega}}(z)\!:=e^{||\mg{\omega}|| z+c^2 ||\bar{\mg{\omega}}||z^{-1}}
\mr{SPA}\!\!\int_0^\infty\!\!\!\!
\frac{e^{-\sum(\omega_i\,t_i+c^2\bar{\omega}_i\,t_i^{-1})}}{(t_r\!-\!t_{r-1})...(t_2\!-\!t_1)(t_1\!-\!z)}dt_1...dt_r
\end{equation}
where {\it SPA} denotes a suitable average of all the $2^{r-1}$ possible integration multipaths that correspond to the $2^{r-1}$ manners in which the variables $t_j$ may circumvent each other on their way from $0$ to  $\infty$. 




%


\section{Multitangents and multizetas.}

The {\it multitangents} and {\it multizetas}, being the transcendental ingredient in the analytical expression of the invariants of identity-tangent diffeos\footnote{and of much else -- they are almost coextensive with the whole field of difference equations.}, deserve a short excursus. But we must begin with a brief reminder about {\it moulds}, which are the proper tool for handling multi-indexed objects of whatever description.

%
\subsection{Mould operations and mould symmetries.}
Moulds are functions of finite sequences  $\mg{\omega}=(\omega_1,...,\omega_r)$ of any length $r\geq 0$,
noted as right-upper indices and rendered, as mute variables, by a plain bold dot $\bu$. Most moulds tend to fall into one or the other of four symmetry classes or types:
\begin{eqnarray*}
M^\bu \mi{symmetral\,(resp. \,alternal)} &\Leftrightarrow& 
\sum_{\mg{\omega} \in \mr{sha}(  \mg{\omega^\prime}, \mg{\omega^{\prime\prime} }) } M^{\mg{\omega}}
= M^{\mg{\omega^\prime}} M^{\mg{\omega^{\prime\prime}}}\;\mi{(resp.\;0)}
\\
M^\bu \mi{symmetrel\,(resp. \,alternel)} &\Leftrightarrow& 
\sum_{\mg{\omega} \in \mr{she}(  \mg{\omega^\prime}, \mg{\omega^{\prime\prime} }) } M^{\mg{\omega}}
= M^{\mg{\omega^\prime}} M^{\mg{\omega^{\prime\prime}}}\;\mi{(resp.\;0)}
\end{eqnarray*}
Here, $\mi{sha}(\mg{\omega^\prime}, \mg{\omega^{\prime\prime} })$
(resp. $\mi{she}(\mg{\omega^\prime}, \mg{\omega^{\prime\prime} })$) 
denotes the set of all sequences $\mg{\omega}$ deducible from $\mg{\omega}^\prime $ and $\mg{\omega}^{\prime\prime}$ under plain (resp. contracting\footnote{
ie allowing for order-compatible, pairwise  contactions 
$(\omega_i^{\prime},\omega_j^{\prime\prime})\mapsto \omega_i^{\prime}+\omega_j^{\prime\prime} $
of elements from the parent sequences.
}  ) shufflings.
Moulds can be {\it multiplied} and {\it composed}\,:
\begin{eqnarray*}
C^\bu=A^\bu\times B^\bu &\Longleftrightarrow & C^{\mg{\omega}}=
\sum_{\mg{\omega^\prime}\mg{\omega^{\prime\prime}}=\mg{\omega}} 
A^{\mg{\omega^{\prime}}}\,B^{\mg{\omega^{\prime\prime}}}
\\
C^\bu=A^\bu\circ B^\bu &\Longleftrightarrow & C^{\mg{\omega}}=
\sum_{\mg{\omega^1}...\mg{\omega^s}= \mg{\omega}} 
A^{\mg{|\omega^1|},..., \mg{|\omega^s|}}\,B^{\mg{\omega^s}}\dots B^{\mg{\omega^s}} \quad\quad
(\mg{|\omega^i|}\not= \mg{\emptyset})
\end{eqnarray*}
with all the predictable relations, including 
$$ (A^\bu\times B^\bu)\circ C^\bu= (A^\bu \circ C^\bu)\times (B^\bu \circ C^\bu)$$.
The multiplication resp. composition unit are the moulds $1^\bu, I^\bu$ defined by:
\begin{eqnarray}    \label{wb1}
\,1^\emptyset:=1\;\; &;&\;\;\;\;\, 1^{\omega_1,\dots, \omega_r}:=0 \quad\;\;\; \mi{if}\;\;\; r\not=0
\\    \label{wb2}
I^{\omega_1}:=1\;\; &;&\;\;\;\;  I^{\omega_1,\dots, \omega_r}:=0 \quad\;\;\;\, \mi{if}\;\;\; r \not= 1
\end{eqnarray}

%
\subsection{Multizetas and multitangents.}
In this subsection, all indices $s_i$ are in $\doN^\ast$ and, to preempt divergence, we (provisionally) assume 
$s_1\not=1$ for multizetas and $s_1,s_r\not=1$ for multitangents.

We first consider two multizeta-valued moulds, the symmetrel $\mi{Ze}^\bu$ and symmetral  $\mi{Za}^\bu$:
\begin{eqnarray}    \label{Ze}
\mr{Ze}^{s_1,...,s_r} &:=& \sum_{n_1>...>n_r>0} n_1^{-s_1}\dots n_r^{-s_r}
\\     \label{Za}
\mr{Za}^{s_1,...,s_r} &:=& \sum_{n_1\geq...\geq n_r>0} n_1^{-s_1}\dots n_r^{-s_r}\, \frac{1}{r_1!}\dots \frac{1}{r_t!}
\end{eqnarray}
In (\ref{Za}), the non-increasing sequences $(n_1,...,n_r)$ involve $t$ clusters of $r_1,...,r_t$ identical integers ($1\leq t \leq r$). Clearly, $ \mi{Za}^\bu=\mi{Ze}^\bu \circ (E^\bu-\mg{1}^\bu)$ with 
\begin{equation}  \label{b5}
(E-1)^{\mg{\emptyset}} =0 \;\;\;\;\;\; \mi{and}\;\;\;\;\;\; 
(E-1)^{s_1,\dots, s_r} = \frac{1}{r!} \;\;\;\; \mi{if}\;\;\;\; 1\leq r\;\;\;\;\; (\forall s_i)
\end{equation}
In a similar vein, we introduce four multitangent-valued moulds, ranging over the four symmetry types:
\[\begin{array}{cccccccccccccccccccccc}
\mr{Te}^\bu & \stackrel{1}{\longrightarrow}  & \mr{Ta}^\bu & \quad\quad
& \mi{symmetrel} & \stackrel{1}{\longrightarrow} & \mi{symmetral}
\\
\downarrow \scriptstyle{3} & & \downarrow \scriptstyle{2} & & \downarrow \scriptstyle{3} & & \downarrow \scriptstyle{2} & &
\\
\mr{Ten}^\bu & \stackrel{4}{\longrightarrow}  & \mr{Tan}^\bu & 
& \mi{alternel} & \stackrel{4}{\longrightarrow} & \mi{alternal}
\end{array}\]
The two upper moulds are defined directly by\footnote{ With the same $r_j$ in (\ref{Ta}) as in (\ref{Za}).}
\begin{eqnarray}    \label{Te}
\mr{Te}^{s_1,...,s_r}(z) &:=& \sum_{n_1>...>n_r} (n_1+z)^{-s_1}\dots (n_r+z)^{-s_r}
\\     \label{Ta}
\mr{Ta}^{s_1,...,s_r}(z) &:=& \sum_{n_1\geq...\geq n_r} (n_1+z)^{-s_1}\dots (n_r+z)^{-s_r} \,\frac{1}{r_1!}\dots \frac{1}{r_t!}
\end{eqnarray}
and the two lower moulds are  their mould logarithms.\footnote{
with the natural definition 
$\mi{logmu}(M^\bu)= M^\bu -\frac{1}{2}M^\bu\times M^\bu +\frac{1}{3}M^\bu\times M^\bu\times M^\bu \dots $.
} Thus:
\begin{eqnarray}  \label{b3}
\mr{Ta}^\bu &=& \mr{Te}^\bu \circ (\mr{E^\bu-1^\bu})
\\  \label{b4}
\mr{Tan}^\bu &=& \mr{Ten}^\bu \circ (\mr{E^\bu-1^\bu})
\\
\label{b1}
\mr{Ten}^\bu &=& \mr{logmu}(\mr{Te}^\bu)
\\  \label{b2}
\mr{Tan}^\bu &=& \mr{logmu}(\mr{Ta}^\bu)
\end{eqnarray}

In the sequel, we shall also consider the inverse $\mi{inv\!Te}^\bu$
of $\mi{Te}^\bu$ (relative to mould multiplication). Clearly
\begin{eqnarray}    \label{invTe}
\mr{invTe}^{s_1,...,s_r}(z) &:=& \sum_{n_1 \leq...\leq n_r} (-1)^r \, (n_1+z)^{-s_1}\dots (n_r+z)^{-s_r}
\end{eqnarray}
All four types of multitangents obviously verify
\begin{equation} \label{pari1}
\mr{T}^{s_1,\dots, s_r}(-z)\equiv (-1)^{s_1+\dots s_r}\,\mr{T}^{s_r,\dots,s_1}(z) \quad\quad 
(\forall \mr{T} \in \{ \mr{Te},\mr{Ta}, \mr{Ten}, \mr{Tan} \})
\end{equation}
In the case of $\mi{Tan}^\bu$, however, due to alternality we have an additional relation
\begin{equation} \label{pari2}
\mr{Tan}^{s_r,\dots,s_1}(z)\equiv (-1)^{r-1}\;\mr{Tan}^{s_1,\dots,s_r}(z)
\end{equation}
Combining (\ref{pari1}), (\ref{pari2}) we get the crucial {\it parity separation} property:
\begin{equation} \label{pari3}
\mr{Tan}^{s_1,\dots,s_r}(-z)\equiv (-1)^{1+\sum d_i}\;\mr{Tan}^{s_1,\dots,s_r}(z)
\quad \quad \mi{with}\; d_i:=s_i-1 
\end{equation}


%
\subsection{Multitangents in terms of multizetas.}
Applying partial fraction decomposition to the series (\ref{Te}), one can easily expand any multitangent $\mi{Te}^\bu$
(hence also any $\mi{Ta}^\bu$, $\mi{Ten}^\bu$ or $\mi{Tan}^\bu$) into a finite sum of utterly elementary monotangents 
${T\!e}^{s_1}$, also known as Eisenstein series. Here is the formula:
\begin{eqnarray}  \label{b6}
\mr{Te}^{s_1,...,s_r} \!\!&=&\!\!\! \sum_{\sigma=2 }^{\sup(s_i)} \mr{Teze}^{s_1,...,s_r}_\sigma \,\mr{Te}^{\sigma}
\;\; = \;\; \sum_{i=1}^{r}\sum_{\sigma_i=2}^{s_i}  \mr{Teze}^{s_1,...,s_r}_{i,\sigma_i} \,\mr{Te}^{\sigma_i}\;\;\; 
\end{eqnarray}
with
\begin{eqnarray*} 
\mr{Teze}^{s_1,...,s_r}_{i,\sigma_i}  \!\!&=&\!\!\! 
\sum_{ \{ { {\sigma_i \leq s_i} \atop  s_j\leq \sigma_j (j\not=i)  } \} }^{\sum \sigma_k = \sum s_k} 
\mr{Ze}^{\sigma_1,...,\sigma_{i-1}}\,
\mr{Ze}^{\sigma_r,...,\sigma_{i+1}}\,
 \prod_{j=1}^{i-1} (-1)^{\sigma_j}\prod_{1\leq j \leq r }^{j\not=i} \frac{(-1)^{s_j}(\sigma_j-1)!}{(\sigma_j-s_j)!(s_j-1)!}
\end{eqnarray*}
or more symmetrically
\begin{eqnarray} \nonumber
\mr{Teze}^{s_1,...,s_r}_{i,\sigma_i}  \!\!&=&\!\!\! 
\sum_{ \{ { {\sigma_i \leq s_i} \atop  s_j\leq \sigma_j (j\not=i)  } \} }^{\sum \sigma_k = \sum s_k} 
\mr{Ze}^{\sigma_1,...,\sigma_{i-1}}\,
(-1)^{s_i-\sigma_i}\,
\mr{viZe}^{\sigma_{i+1},...,\sigma_{r}}\,
 \prod_{1\leq j \leq r }^{j\not=i} \frac{(\sigma_j-1)!}{(\sigma_j-s_j)!(s_j-1)!}
 \\ \label{b6bis}
 \mr{viZe}^{s_1,...,s_r} \!\!&=&\!\!\! (-1)^{s_1+...s_r}\,\mr{Ze}^{s_r,...,s_1}
\end{eqnarray}
The leading monotangent $\mi{Te}^{1}(z)=\frac{\pi}{\mi{tan}(\pi z )}$ generates all others under differentiation, and admits the following northern and southern expansions:
\begin{eqnarray}   \label{b7}
\mr{Te}_{\mr{no}}^1(z) &=& -\pi i -2 \pi i\sum_{0<n} e^{+2 \pi i\,n\,z}  \hspace{6. ex} \mi{if}\;\; \Im (z) >0
\\    \label{b7b}
\mr{Te}_{\mr{so}}^1(z) &=& +\pi i +2 \pi i\sum_{0<n} e^{-2 \pi i\,n\,z}  \hspace{6. ex} \mi{if}\;\; \Im (z) <0
\end{eqnarray}
Since $\mr{Te}^{s_1}(z)= \frac{(-1)^{s_1-1}}{(s_1-1)!}\,\partial_z^{s_1-1}\, \mr{Te}^1(z) $, this yields

\begin{eqnarray}   \label{b8}
\mr{Te}^{s_1}(z) &=& \sum_{\omega\in \Omega^{\mp}} \, \mr{Te}^{s_1}_\omega\, e^{-\omega z}  
\hspace{6. ex} \mi{on\, each\, half\textup{-}plane } \;\; \pm \Im(z) >0
\end{eqnarray}
with
\begin{eqnarray}   \label{b9}
\mr{Te}^{s_1}_\omega &=&  \mr{sign}(\Im (\omega))\, 2\pi i\,\frac{\omega^{s_1-1}}{(s_1-1)!}
\quad \mi{and} \quad \Omega^\mp =  2\pi i \doZ^{\mp}
\end{eqnarray}
All the above amounts to a simple procedure for calculating the Fourier expansions, north and south,
of the four classes of multitangents. The classes 
$\mi{T\!e}^\bu$ and $\mi{T\!an}^\bu$ shall be of direct concern to us:
\begin{eqnarray}      \label{B1}
\mr{Te}^\bu_{\mr{no}}(z) = \sum_{\omega\in \Omega^{-}} \mr{Te}^\bu_{\omega}\;\, e^{-\omega\,z}\;\;\;
&\quad;\quad &
\;\, \mr{Te}^\bu_{\mr{so}}(z) = \sum_{\omega\in \Omega^{+}} \mr{Te}^\bu_{\omega}\;\, e^{-\omega\,z}
\quad
\\     \label{B2}
\mr{Tan}^\bu_{\mr{no}}(z) = \sum_{\omega\in \Omega^{-}} \mr{Tan}^\bu_{\omega}\;\, e^{-\omega\,z}
&\quad;\quad &
\mr{Tan}^\bu_{\mr{so}}(z) = \sum_{\omega\in \Omega^{+}} \mr{Tan}^\bu_{\omega}\;\, e^{-\omega\,z}
\quad
\end{eqnarray}
{\bf Remark: advantages of  $\mi{Tan}^\bu $  over $\mi{Te}^\bu $.}
\\
Unlike the symmetrel multitangents $\mi{Te}^\bu$, their alternal counterparts $\mi{Tan}^\bu$ admit of no simple, direct expansions of type (\ref{Te}), and their expression as superpositions of $\mi{Te}^\bu$ is very involved, as shown by the formulae of \S5.1. The picture changes, however, after reduction into monotangents and symmetrel linearisation of the resulting multizetas: it is now  $\mi{Tan}^\bu$ that gives rise, by and large, to the simpler expansions, as shown by the Table at the beginning of \S5.2.  Then   $\mi{Tan}^\bu$ possesses a second advantage: that of having a definite parity, which depends only on the total degree $\sum d_i$: cf (\ref{pari3}) above. Thirdly, we shall see in the sections \S3 and \S4 that, when it comes to expressing the collectors or the invariants, $\mi{Tan}^\bu$ leads to decidedly simpler formulae than 
$\mi{Te}^\bu$, as immediately apparent from a comparison of (\ref{vivi10}) with  (\ref{vivi4})-(\ref{vivi5})  or of (\ref{catur4}) with 
(\ref{catur3}).


%
\subsection{Multitangents in terms of resurgent monomials.}
There exists an alternative, {\it resurgent} approach to multitangent reduction. In the convergent (i.e. $s_1,s_r\not=1$) and non-ramified (i.e. $s_j\in \doN^\ast$ rather than $\doQ^\ast$) case, it hardly improves on the above procedure  (see \S2.3) but in the general case, especially when we go over to fractional indices $s_j$, the resurgent approach becomes the more flexible of the two methods and even, in a sense, the only practical one. For clarity, though, we first sketch the alternative method under retention of the two simplifying assumptions: no divergence, no ramification.

We begin by constructing two resurgent-valued, symmetrel, mutually inverse moulds, first in the {\it formal model}, via the induction:

\begin{eqnarray}  \label{bab1}
\mr{\tilde{S}e}^\bu\!(z) &=& 
\frac{e^{\partial_z}}{(1-e^{\partial_z})}\, 
\Big(\mr{\tilde{S}e}^\bu\!(z) \times \mr{\stackrel{}{Je}}^\bu\!(z)\Big)
\\  \label{bab2}
\mr{inv\tilde{S}e}^\bu\!(z) &=& 
\frac{e^{\partial_z}}{(e^{\partial_z}-1)}\, 
\Big(\mr{\stackrel{}{Je}}^\bu\!(z)\times\mr{inv\tilde{S}e}^\bu\!(z)\Big)
\end{eqnarray}
 with the elementary mould $\mi{\stackrel{}{J\!e}\!}^{\bu}(z)$ defined by
\begin{equation}    \label{wb5}
\mr{Je}^{\emptyset}(z):= 0\;\;\; ;\;\;\;
\mr{Je}^{s_1}(z):= z^{-s_1}\;\;; \;\;
\mr{Je}^{s_1,...,s_r}(z):= 0\;\;\; (\forall\, r\geq2)
\end{equation}
Together with the conditions 
\begin{equation}    \label{wb6}
\mr{\tilde{S}e}^\emptyset\!(z)=\mr{inv\tilde{S}e}^\emptyset\!(z)=1\,, \,
\mr{\tilde{S}e}^{s_1,...,s_r}\!(\infty)=\mr{inv\tilde{S}e}^{s_1,...,s_r}\!(\infty)=0 \;\;(\forall r\geq1)\;\;\;
\end{equation}
the induction (\ref{bab1})-(\ref{bab2}) uniquely defines each $\mi{\tilde{S}\!e}^{\mg{s}}\!(z)$ 
and each $\mi{inv\!\tilde{S}\!e}^{\mg{s}}\!(z)$ as a constant-free power series in $z^{-1}$. 

In the {\it convolutive model} the induction becomes
\begin{eqnarray}  \label{bab3}
\mr{\hat{S}e}^{s_1,...,s_r}(\zeta) \!&=&\! \frac{e^{-\zeta}}{(1-e^{-\zeta})}\,
\int_0^{\zeta} \mr{\hat{S}e}^{s_1,...,s_{r-1}}\!(\zeta\!-\!\zeta_r) \;\,
\frac{\zeta_r^{s_r-1}}{\Gamma(s_r)}\; d\zeta_r
\\  \label{bab4}
\mr{inv\hat{S}e}^{s_1,...,s_r}(\zeta) \!&=&\! \frac{e^{-\zeta}}{(e^{-\zeta}-1)}\,
\int_0^{\zeta} 
\frac{\zeta_1^{s_1-1}}{\Gamma(s_1)}\;\,
\mr{inv\hat{S}e}^{s_2,...,s_{r}}\!(\zeta\!-\!\zeta_1) \; d\zeta_1
\end{eqnarray}

Lastly, in the {\it geometric models} $+$ and $-$ ({\it east} and {\it west}), corresponding to Laplace integration along the axes
$\mi{arg}(\zeta)=0$ and $\mi{arg}(\zeta)=\pi$, we get
\begin{eqnarray}  \label{bb1}
\mr{Se}_{+}^{s_1,...,s_r}(z) &:=&\sum_{0<n_r<...<n_1} \hspace{5.ex} (n_1+z)^{-s_1}\dots (n_r+z)^{-s_r}
\\  \label{bb2}
\mr{Se}_{-}^{s_1,...,s_r}(z) &:=&\sum_{n_1\leq..\leq n_r\leq 0} (-1)^r \,(n_1+z)^{-s_1}\dots (n_r+z)^{-s_r}
\\  \label{bb3}
\mr{invSe}_{+}^{s_1,...,s_r}(z) &:=&\sum_{0 < n_1\leq...\leq n_r} (-1)^r \,(n_1+z)^{-s_1}\dots (n_r+z)^{-s_r}
\\  \label{bb4}
\mr{invSe}_{-}^{s_1,...,s_r}(z) &:=&\sum_{n_r<..< n_1 \leq 0} \hspace{7.ex}   (n_1+z)^{-s_1}\dots (n_r+z)^{-s_r}
\end{eqnarray}

From the structure of the induction (\ref{bab1}), one directly (without calculation) infers that the $\mi{ \tilde{S}\!e}^\bu(z)$ must verify resurgence equations of the form\footnote{we drop the tilde for simplicity.}
\begin{eqnarray}  \label{bb5}
+2\pi i\,\Delta_\omega\; \mr{Se}^\bu(z)  &=& \mr{Ten}^\bu_\omega \times \mr{Se}^\bu(z)     
 \hspace{7. ex}  ( \forall \omega \in \Omega^+ = 2\pi i \doZ^+)
\\  \label{bb6}
-2\pi i\,\Delta_\omega\; \mr{Se}^\bu(z)   &=& \mr{Ten}^\bu_\omega \times \mr{Se}^\bu(z)     
 \hspace{7. ex}  ( \forall \omega \in \Omega^- = 2\pi i \doZ^-)
\\  \label{bb7}
\Delta^{+}_\omega\; \mr{Se}^\bu(z)  &=& \mr{Te}^\bu_\omega \times \mr{Se}^\bu(z) 
  \hspace{8.ex }  ( \forall \omega \in  \Omega^+ = 2\pi i \doZ^{+})
\\  \label{bb8}
\Delta^{-}_\omega\; \mr{Se}^\bu(z)  &=& \mr{Te}^\bu_\omega \times \mr{Se}^\bu(z)  
 \hspace{8.ex}  ( \forall \omega \in   \Omega^- = 2\pi i \doZ^{-})
\end{eqnarray}
with scalar-valued moulds $\mi{Ten}_\omega^\bu$ (alternel) and  $\mi{Te}_\omega^\bu$ which, for the moment, need not  bear any relation to their namesakes in (\ref{B1}).
However, applying the axis-crossing identity (\ref{cross+}) to (\ref{bb7}) with $\theta=+\frac{\pi}{2}$ and the reverse identity
(\ref{cross-}) to (\ref{bb8}) with $\theta=-\frac{\pi}{2}$ , and minding the fact that
$$
\mr{Se}^\bu_{\frac{\pi}{2}\pm\epsilon}=\mr{Se}^\bu_{\mp}\;\;\; \mi(inversion!)
\quad ;\quad
\mr{Se}^\bu_{-\frac{\pi}{2}\pm\epsilon}=\mr{Se}^\bu_{\pm}\;\;\; \mi(no\;\,inversion!)
$$
we find respectively
\begin{eqnarray}  \label{bb10}
 \mr{Te}_{\mr{so}}^\bu(z)\times  \mr{Se}_{-,\mr{so}}^\bu(z) =  \mr{Se}_{+,\mr{so}}^\bu(z)
 & \mi{with} &
  \mr{Te}_{\mr{so}}^\bu(z) =  \sum_{\omega\in \Omega^+} \mr{Te}^\bu_\omega\; e^{-\omega z}
\\   \label{bb11}
\mi{Te}_{\mr{no}}^\bu(z)\times  \mr{Se}_{-,\mr{no}}^\bu(z) = \mr{Se}_{+,\mr{no}}^\bu(z)
 & \mi{with} &
 \mr{Te}_{\mr{no}}^\bu(z) = \sum_{\omega\in \Omega^-} \mr{Te}^\bu_\omega\; e^{-\omega z}
\end{eqnarray}

Thus, whether ``north" or ``south", we arrive at the elementary identity
\begin{eqnarray}\label{bb12}
 \mr{Te}^\bu(z) &=&  \mr{Se}_{+}^\bu(z) \times \mr{invSe}_{-}^\bu(z)
\end{eqnarray}
which of course can also be directly derived from the definitions (\ref{Te}), (\ref{bb1}), (\ref{bb4}). But we get an interesting extra -- namely, that the moulds $\mi{Te}_\omega^\bu$ of  (\ref{bb7}) and  (\ref{bb8}) coincide with those of (\ref{B1}) in the preceding subsection. If we now interpret the resurgence equations (\ref{bb5})-(\ref{bb8}) in the convolutive model, we get an alternative expression of  $\mi{Te}_\omega^\bu$ and   $\mi{Ten}_\omega^\bu$
(and hence $\mi{Ta}_\omega^\bu$,  $\mi{Tan}_\omega^\bu$) as finite integrals in the $zeta$-plane, which translates, after some work, into fast-convergent power series. This will stand us in good stead in the {\it divergent} and above all in the {\it ramified} cases. But we must first devote a short aside to the question of parity.
 

%
\subsection{Respecting and harnessing parity.}
When it comes to calculating what will turn out to be most basic and useful of all four multitangents, namely
$\mi{Tan}^\bu$, all the above procedures must be adjudged {\it wasteful} in the sense that they derive the parity-separating
 $\mi{Tan}^\bu$ from the parity-mixing $\mi{Te}^\bu$.
 To remedy this, we shall
 replace the nearly odd function $\frac{e^\partial}{1-e^\partial} $ of (\ref{bab1})-(\ref{bab2}) by the exactly odd $H(\partial)$:
\begin{equation}\label{baab0}
H(\partial):= \frac{e^\partial}{1-e^\partial} +\frac{1}{2} =
\frac{1}{2}\, \frac{1+e^\partial}{1-e^\partial}= -\frac{1}{2}\; \mr{cotan}(\frac{\partial}{2})
\end{equation}
and define two mutually inverse (but no longer exactly symmetrel) moulds 
$\mi{\tilde{S}\!ee}^\bu$ and  $\mi{inv\!\tilde{S}\!ee}^\bu$ by the tweaked induction
\begin{eqnarray}  \label{baab1}
\mr{\tilde{S}ee}^{\bu}\!(z) &=& +H(\partial_z)\, \Big(\mr{\tilde{S}ee}^{\bu}\!(z) \times {\stackrel{}{\mr{Je}}}^{\bu}(z)\Big)
\\  \label{baab2}
\mr{inv\tilde{S}ee}^\bu\!(z) &=& -H(\partial_z)\,\, \Big(\mr{\stackrel{}{Je}}^\bu\!(z)\times\mr{inv\tilde{S}ee}^\bu\!(z)\Big)
\end{eqnarray}
Clearly, $\mi{\tilde{S}\!ee}^{s_1,...,s_r}\!(z)$ and  $\mi{inv\!\tilde{S}\!ee}^{s_1,...,s_r}\!(z)$ are {\it even} (resp. {\it odd}) power series in $z^{-1}$ iff the total ``degree'' $(\sum s_i)\m r$ is  {\it even} (resp. {\it odd}).

As in the preceding subsection, the induction immediately leads to resurgence equations of the form:
\begin{eqnarray}  \label{baab5}
+2\pi i\,\Delta_\omega\; \mr{See}^\bu(z)  &=& \mr{Teen}^\bu_\omega \times \mr{See}^\bu(z)     
 \hspace{7. ex}  ( \forall \omega \in \Omega^+ = 2\pi i \doZ^+)
\\  \label{baab6}
-2\pi i\,\Delta_\omega\; \mr{See}^\bu(z)   &=& \mr{Teen}^\bu_\omega \times \mr{See}^\bu(z)     
 \hspace{7. ex}  ( \forall \omega \in \Omega^- = 2\pi i \doZ^-)
\\  \label{baab7}
\Delta^{+}_\omega\; \mr{See}^\bu(z)  &=& \mr{Tee}^\bu_\omega \times \mr{See}^\bu(z) 
  \hspace{8.ex }  ( \forall \omega \in  \Omega^+ = 2\pi i \doZ^{+})
\\  \label{baab8}
\Delta^{-}_\omega\; \mr{See}^\bu(z)  &=& \mr{Tee}^\bu_\omega \times \mr{See}^\bu(z)  
 \hspace{8.ex}  ( \forall \omega \in   \Omega^- = 2\pi i \doZ^{-})
\end{eqnarray}
Moreover, (\ref{baab0}) implies that the moulds $\mi{\tilde{S}\!e}^\bu\!(z)$ and  $\mi{\tilde{S}\!ee}^\bu\!(z)$ are related as follows
\begin{eqnarray} \label{baab9}
\mr{\tilde{S}e}^\bu&=&\Big(\mr{\tilde{S}ee}^\bu \circ \mr{\stackrel{}{D}}_{-\frac{1}{2}}^\bu  \Big)\times
\Big(\mg{1}^\bu -\frac{1}{2} \mr{\stackrel{}{Je}}^\bu\!(z)\circ \mr{\stackrel{}{D}}_{-\frac{1}{2}}^\bu  \Big)
\\   \label{baab10}
\mr{\tilde{S}ee}^\bu&=&\Big(\,\mr{\tilde{S}e}^\bu\, \circ \mr{\stackrel{}{D}}_{+\frac{1}{2}}^\bu  \Big)\times
\Big(\mg{1}^\bu +\frac{1}{2} \mr{\stackrel{}{Je}}^\bu\!(z)\circ \mr{\stackrel{}{D}}_{+\frac{1}{2}}^\bu  \Big)
\end{eqnarray}
with
\begin{equation}  \label{baab11}
\mr{\stackrel{}{D}}_a^\emptyset :=0\quad;\quad \mr{\stackrel{}{D}}_a^{s_1,\dots,s_r}:=a^{r-1}\quad (\forall r\geq 1)
\end{equation}
However, the elementary second factors on the right-hand sides of (\ref{baab9}) and  (\ref{baab9}) are convergent rather than resurgent-valued. As a consequence, the corresponding right and left germs coincide, and if we set
\begin{eqnarray}\label{baab12}
 \mr{Tee}^\bu(z) &=&  \mr{See}_{+}^\bu(z) \times \mr{invSee}_{-}^\bu(z)
\end{eqnarray}
we shall have 
\begin{eqnarray*}
 \mr{Te}^\bu(z) &=& \mr{Se}_{+}^\bu(z) \times \mr{invSe}_{-}^\bu(z) \;=\;
  \Big(\mr{See}_{+}^\bu(z) \times \mr{invSee}_{-}^\bu(z)\Big)\circ \mr{\stackrel{}{D}}_{-\frac{1}{2}}^\bu
\end{eqnarray*}
and therefore
\begin{eqnarray}\label{baab13}
 \mr{Te}^\bu(z) &=&   \mr{Tee}^\bu(z) \circ  \mr{\stackrel{}{D}}_{-\frac{1}{2}}^\bu
\end{eqnarray}
Postcomposing this by $\mg{E}^\bu-\mg{1}^\bu$ and defining an elementary, purely length-dependent mould $K^\bu$ by:
\begin{equation} \label{baab14}
\sum_{0\leq r} K^{s_1,\dots, s_r}\; t^r := 2\,\tan(\frac{t}{2})\quad\quad\quad\quad (\forall s_i)
\end{equation}
we get:
\begin{eqnarray}\label{baab15}
 \mr{Ta}^\bu(z) &=&   \mr{Tee}^\bu(z) \circ K^\bu
\end{eqnarray}
Lastly, taking in the above identity the mould logarithm of both sides, we arrive at:
\begin{eqnarray}\label{baab16}
 \mr{Tan}^\bu(z) &=&   \mr{Teen}^\bu(z) \circ K^\bu
 \\ \label{baab17}
  \mr{Tan}_\omega^\bu(z) &=&   \mr{Teen}_\omega^\bu(z) \circ K^\bu 
  \quad \quad \forall \omega\in\Omega=2\,\pi i \,\doZ^\ast)
\end{eqnarray}
We have thus expressed the basic alternal multitangents $\mi{Tan}^\bu$ directly in terms of  auxiliary multitangents 
 $\mi{Teen}^\bu$ which, though belonging to none of the basic symmetry types, 
 share with  $\mi{Tan}^\bu$  the same crucial property of parity separation:
 as functions of $z$, both $\mi{Teen}^{s_1,..,s_r}\!(z)$ and $\mi{Tan}^{s_1,..,s_r}\!(z)$  have the same parity as 
 that of the number $1-r+\sum s_i$.


%
\subsection{The divergent case. Normalisation.}
If we now drop the condition that ensured convergence, namely $s_1, s_r\not=1$, and yet insist on retaining all properties and symmetries of our moulds, we must do two things to our infinite series: {\it truncate} them and {\it correct} them. Concretely, we must set
\[\begin{array}   {cclllllcc}
\mr{Te}^\bu\!(z)
\!\!\!\!&\!:=\!&\!\!
\lim_{k\rightarrow \infty} \mr{Te}^\bu_k(z) 
\!\!\!\!&\!:=\!&\!\!
 \lim_{k\rightarrow \infty}\; \mr{invcoSe}^\bu_k \times \mr{doTe}^\bu_k(z) \times  \mr{coSe}^\bu_k
 \\ [1. ex]
\mr{Se}_{\pm}^\bu\!(z)
\!\!\!\!&\!:=\!&\!\!
\lim_{k\rightarrow \infty}  \mr{Se}_{k,\pm}^\bu\!(z) 
\!\!\!\!&\!:=\!&\!\!
 \lim_{k\rightarrow \infty}\; \mr{invcoSe}_{k}^\bu\times \mr{doSe}_{k,\pm}^\bu\!(z)
\\    [1. ex]
\mr{invSe}_{\pm}^\bu\!(z)
\!\!\!\!&\!:=\!&\!\!
\lim_{k\rightarrow \infty}  \mr{invSe}_{k,\pm}^\bu\!(z) 
\!\!\!\!&\!:=\!&\!\!   
 \lim_{k\rightarrow \infty}\; \mr{invdoSe}_{k,\pm}^\bu\!(z)  \times \mr{coSe}_{k}^\bu 
\end{array}\]
Here, the symmetrel {\it dominant} factors $\mi{Te}^\bu$, $\mi{doSe}_{k,\pm}^\bu$,  
$\mi{invdoSe}_{k,\pm}^\bu$ are defined 
as in (\ref{Te}) and (\ref{bb1})-(\ref{bb4}) but with sums truncated at $\pm k$ instead of  $\pm \infty$. Thus
\begin{eqnarray}    \label{diver2}
\mr{doTe}_k^{s_1,...,s_r}(z) &:=& \sum_{-k\leq n_r<...<n_1\leq k} (n_r+z)^{-s_r}\dots (n_1+z)^{-s_1}
\;\; (\forall s_i)
\end{eqnarray}
As for the symmetrel, $z$-constant {\it corrective} factors $\mi{coSe}_{k\pm}^\bu$ and  $\mi{invcoSe}_{k\pm}^\bu$, 
their definition reduces to
\begin{eqnarray} \label{diver3}
\mr{coSe}_k^{s_1,...,s_r}&:=& \frac{(\; c \,+\log k\,)^r}{r!} 
\,\quad\quad \mi{if}\;\; (s_1,...,s_r)=(1,...,1)
\\    \label{diver4}
\mr{invcoSe}_k^{s_1,...,s_r}&:=& \frac{(-c-\log k )^r}{r!} 
\,\quad\quad \mi{if}\;\; (s_1,...,s_r)=(1,...,1)
\\ [1.ex]   \label{diver5}
\mr{coSe}_k^{s_1,...,s_r}&=&\!\!\mr{invcoSe}_k^{s_1,...,s_r} :=\;0 
\quad \mi{if}\;\; (s_1,...,s_r)\not=(1,...,1)
\end{eqnarray}

In the formal model, the resurgent-valued moulds  $\mi{\tilde{S}\!e}^{\bu}$ and $\mi{inv\!\tilde{S}\!e}^{\bu}$ are still 
uniquely defined by the induction (\ref{bab1})-(\ref{bab2}) together with the condition 
\begin{equation} \label{diver6}
\mr{\tilde{S}e}^{\mg{s}}(z)\;\; ,\;\; \mr{inv\tilde{S}e}^{\mg{s}}(z)\;\;\in 
\doQ[[z^{-1}]]\otimes \doQ[(c+\log z)] \stackrel{.}{-} \doQ
\quad\quad (\forall \mg{s}\not= \mg{\emptyset})
\end{equation}
The normalising condition, in other words, is that $\mr{\tilde{S}e}^{\mg{s}}(z)$ and $\mr{inv\tilde{S}e}^{\mg{s}}(z)$, 
as formal series in $z^{-1}$ and  polynomials in the bloc $(\mi{c\p log\,z })$, should have no constant term.

In the sectorial models, the $c$-normalisation implies:
\begin{equation}    \label{wb20}
{Se}_{\pm}^{\stackrel{r\,\mr{times}}{\overbrace{1,...,1}}}(0)=\frac{(\gamma-c)^r}{r!}
\quad ; \quad
{invSe}_{\pm}^{\stackrel{r\,\mr{times}}{\overbrace{1,...,1}}}(0)=\frac{(c-\gamma)^r}{r!}
\end{equation}
with 
\begin{equation}    \label{wb21}
\gamma=\lim_{k\rightarrow \infty}(1+\frac{1}{2}+...+\frac{1}{k}-\log k)=0.577215...=\mi{Euler\; constant}
\end{equation}

For multitangents, we may still formally apply the procedure 
(\ref{b6})-(\ref{b6bis}) of \S2.3 to reduce them into combinations of monotangents and mutizetas, but this
time we are liable to get formally divergent multizetas. 
The  $c$-normalisation then amounts to setting $\zeta(1)=\mi{Ze}^1:=\gamma-c$ and to adopting for all divergent multizetas\footnote{i.e. for all multizetas with initial index $s_1=1$.} the unique {\it symmetrel extension} compatible with that initial choice.\footnote{
Thus $\mi{Ze}^{1,1}:=-\frac{1}{2}\mi{Ze}^2+\frac{1}{2}(\gamma-c)^2 $\,,\,
$\mi{Ze}^{1,2}:=-\mi{Ze}^{2,1}-\mi{Ze}^3+(\gamma-c)\,\mi{Ze}^2$\, etc. There exist simple formulae for calculating the symmetrel extension of all multizetas relative to any given choice of $\mi{Ze}^1$.
}

There are two natural choices for the normalisation constant $c$\,: 
\\
(i) Either we set  $c=0$, in which case we eschew $\gamma$ in the formal model but at the cost of introducing it in the convolutive and sectorial models. It also complicates the definition of the multitangents and multizetas, since it forces us to set $\mi{Ze}^1=\gamma$, which however is not entirely unnatural, in view of the formula
\begin{equation} \label{diver7}
\sigma\,\Gamma(\sigma)= \exp\Big(-\gamma\,\sigma+\sum_{2\leq n} (-1)^n\frac{\zeta(n)}{n} \sigma^n\Big)
\end{equation}
(ii) Or we set $c=\gamma$, which forcibly introduces $\gamma$ into the formal model but rids us of it everywhere else, including in the definition of multitangents and multizetas, since it amounts to setting $\mi{Ze}^1=0$. This shall be our preferred choice.


%
\subsection{The ramified case.}
This is the case of fractional {\it weights} $s_i$  in ${p}^{-1}\doN^\ast$ and no longer in $\doN^\ast$.
Everything carries over to that case, except the {\it finite} reduction of multitangents into monotangents and multizetas.

The formulae (\ref{b6})-(\ref{b6bis}) still make formal sense but lead to expansions which are not only infinite but also divergent. When properly re-summed, they yield the correct expressions, but from the point of view of calculational efficiency, this approach is worthless.

Of course, straightforward Fourier analysis in the upper and lower halves of the $z$-plane would yield the coefficients 
$\mi{Te}^\bu_\omega$ along with all the others, but not in the form of nice convergent series, and again at great cost.

The resurgence approach of \S2.4 and \S2-5, on the other hand, survives ramification without any modification. When pursued to the end, this approach even leads to some sort of functional equation for multizetas, that is to say, to something resembling the classical relation between $\zeta(s)$ and $\zeta(1-s)$.





%


%
\section{Collectors and connectors in terms of $f$.}



%
\subsection{Operator relations.}

To the composition identity between germs 
$ f^\ast=l^{-1} \circ  f^\ast\circ f \equiv l^{-1} \circ f^\ast\circ l\circ g $ there answers the operatorial relation\footnote{
To diffeos $f, g ...$ we associate the operators $F, G ...$ of postcomposition by $f,g...$}
\begin{equation} \label{iter}
F^\ast\;=\;G\,\, F_{:1}^\ast \quad \quad \mi{with} \quad F_{:1}^\ast:=L\,F^\ast L^{-1} \quad\quad
\end{equation}
To solve (\ref{iter}) while respecting the basic symmetry between $f,g$ and $f^{-1},g^{-1}$, we set as our basic `infinitesimals' the following operators
\begin{eqnarray}\label{iter0}
\fmi{G}_{:n}^+ &:=& L^n.(G\;-\;1).L^{-n} \hspace{7.ex} (n_i\in \doZ)
\\[0.7 ex]   \label{iter00}
\fmi{G}_{:n}^- &:=& L^n.(G^{-1}\!\!-1).L^{-n} \hspace{7.ex} (n_i\in \doZ)
\end{eqnarray}
This leads straightaway to the formal expansions
\begin{eqnarray}\label{iter1}
F^\ast_{\pp}&=&1+\sum_{1\leq r}\quad\sum_{0\leq n_r<...<n_1}\fmi{G}^{\pp}_{:n_r}\dots \fmi{G}^{\pp}_{:n_1}
\hspace{7.ex} (n_i\in \doZ)
\\   \label{iter2}
F^\ast_{\mm}&=&1+\sum_{1\leq r}\quad\sum_{n_1<...<n_r<0}\fmi{G}^{\mm}_{:n_r}\dots \fmi{G}^{\mm}_{:n_1}
\hspace{7.ex} (n_i\in \doZ)
\\[1.ex]   \label{iter3}
{^\ast\!}{F}_{\pp}&=&1+\sum_{1\leq r}\quad\sum_{0\leq n_1<...<n_r}\fmi{G}^{\mm}_{: n_r}\dots \fmi{G}^{\mm}_{: n_1}
\hspace{7.ex} (n_i\in \doZ)
\\   \label{iter4}
{^\ast\!}{F}_{\mm}&=&1+\sum_{1\leq r}\quad\sum_{n_r<...<n_1<0}\fmi{G}^{\pp}_{:n_r}\dots \fmi{G}^{\pp}_{:n_1}
\hspace{7.ex} (n_i\in \doZ)
\end{eqnarray}
\begin{eqnarray}\label{iter5}
\mg{\Pi}^+\,:={^\ast\!}{F}_{\mm}.F^\ast_{\pp}&=&1
+\sum_{1\leq r}\quad\sum_{n_r<...<n_1}\fmi{G}^{\pp}_{:n_r}\dots \fmi{G}^{\pp}_{:n_1}
\hspace{4.ex} (n_i\in \doZ)
\\   \label{iter6}
\mg{\Pi}^{-}\,:={^\ast\!}{F}_{\pp}.F^\ast_{\mm}&=&
1+\sum_{1\leq r}\quad\sum_{n_1<...<n_r}\fmi{G}^{\mm}_{:n_r}\dots \fmi{G}^{\mm}_{:n_1}
\hspace{4.ex} (n_i\in \doZ)
\end{eqnarray}
For standard diffeos $f$, the above expansions for $F^\ast, ^{\ast\!}F$ (resp. $\Pi^{\pm1}$ ) are easily shown to {\it converge} when they are made to act on test functions that are defined on suitably extended U-shaped domains (resp. on suitably distant  half-planes $|\Im(z)|\gg 1$).
 
The challenge is now to extract from these expansions (- first in the standard, then in the general case -) {\it theoretically appealing}, {\it analytically transparent}, and {\it computationally manageable} expressions for (in that order) the {\it collectors}, {\it connectors}, and {\it invariants}.


%
\subsection{The direct, non-symmetrical scheme.}
To express our `infinitesimals' $\fmi{G}^{\pm}_{:n}$ in terms of the diffeo's Taylor coefficients, we first set
\begin{equation}                        \label{vivi1}
\fmi{G}^{\pm}_{:n} = \sum_{1\leq k}\; \frac{1}{k!}\; \big( g^{\pm1}(z+n) - (z+n) \big)^k\; \partial_z^k
\end{equation}
and then
\begin{eqnarray} \!\!\!\!\!                        \label{vivi2}
\big(g(z)-z\big)^n= \sum_{1\leq d} g^{+}_{n,1+d}\,z^{-d} \!&=&\! \sum_{2\leq s} g^{+}_{n,s}\,z^{-s+1} \quad\quad 
(d\!=\!\mi{``degree"})
\\   \!\!\!\!\!                        \label{vivi3}
\big(g^{-1}\!(z)-z\big)^n= \sum_{1\leq d} g^{-}_{n,1+d}\,z^{-d} \!&=&\! \sum_{2\leq s} g^{-}_{n,s}\,z^{-s+1} \quad\quad 
(s\!=\!\mi{``weight"})
\end{eqnarray}
Next, to account for the action of the {\it derivation} operators $\partial_z$ implicit in the definition of the {\it substitution} operators
$ \fmi{G}^{\pm}_{:n}$, we require integers $\delta^\bu_\bu$ defined by 
\begin{eqnarray}                        \label{vivi333}
\sum_{\sum(l_i-n_i)=1} \!\!\! \delta^{\,l_1,...,\,l_r}_{n_1,...,n_r}
\, x_1^{l_1} \dots x_r^{l_r} &\equiv& x_1^{n_2}.(x_1\!+\!x_2)^{n_3}\dots (x_1\!+\!\dots x_{r-1})^{n_r}
\end{eqnarray}
Letting the operators on both sides of (\ref{iter5}) resp. (\ref{iter6}) act on the test function $z$, and harvesting all $r$-linear summands, we find the sought-after expansions for the collectors $\mg{\gop}^{\pm}$:
\begin{eqnarray}                         \label{vivi4}
\!\!\! \!\!\!\!\!\! \!\!\!
\mg{\gop}(z)  \!\!&\!\!=\!\!&\!\!
z\!+\!\sum_{1\leq r}\!\!\!
\sum_{  {0 \leq l_i \atop 1\leq n_i }}^{n_i+l_i\leq s_i}
\!\!\!\!(-1)^{n\!-\!1}
\! \, \delta^{\,l_1,...,\,l_r}_{n_1,...,n_r}\, \mathbf{Te}^{s_1,..., s_r}(z)
\!\!  \prod_{1\leq i \leq r} \!\!
\frac{(s_i\!-\!1)!\,g^{+}_{n_i,s_i\!-\!l_i\!+\!1}}{(s_i\!-\!l_i\!-\!1)!} 
  \\ [1.5 ex]                        \label{vivi5}
\!\!\!\!\!\!\!\!\! \!\!\!
\mg{\gop}^{-1}(z)  \!\!&\!\!=\!\!&\!\!
z\!+\!\sum_{1\leq r}\!\!\!
\sum_{  {0 \leq l_i \atop 1\leq n_i }}^{n_i+l_i\leq s_i}
\!\!\!\!(-1)^{n\!-\!1}
\! \, \delta^{\,l_1,...,\,l_r}_{n_1,...,n_r}\, \mathbf{Te}^{s_1,..., s_r}(z)
\!\!  \prod_{1\leq i \leq r} \!\!
\frac{(s_i\!-\!1)!\,g^{-}_{n_i,s_i\!-\!l_i\!+\!1}}{(s_i\!-\!l_i\!-\!1)!}   
\end{eqnarray}
with $n:=n_1+...n_r$.


\subsection{The indirect, symmetrical scheme.}
As already pointed out in \S1.4, the only way to fully respect the underlying symmetries between $g,\mg{\gop}$ and 
 $g^{-1},\mg{\gop}^{-1}$ is to switch to the formal generators $g_\ast,\mg{\gop}_\ast$. That means expressing the substitution operators  $\fmi{G}_{:n}$ in terms of the (generically divergent) series
\begin{equation}                        \label{vivi7}
g_\ast(z)= \sum_{1\leq d} g_{\ast 1+d}\,z^{-d} = \sum_{2\leq s} g_{\ast s}\,z^{-s+1} \quad\quad 
(d\!=\!\mi{degree},s\!=\!\mi{weight})
\end{equation}
and setting
\begin{equation}                        \label{vivi6}
\fmi{G}_{:n} = \sum_{1\leq k}\; \frac{1}{k!}\; \big( g_\ast(z+n)\; \partial_z \big)^k
\end{equation}

This time around, we require much simpler coefficients $\delta^\bu$ and $\delta_1^\bu$ with only one sequence of (upper) indices:
\begin{eqnarray}                        \label{vivi8}
\sum_{l_i\geq 0\,,\, \sum l_i=r-1}  \delta^{l_1,...,l_r}\, x_1^{l_1} \dots x_r^{l_r} &\equiv& x_1.(x_1+x_2)\dots (x_1+\dots+x_{r-1})
\\                        \label{vivi9}
\sum_{l_i\geq 0\,,\, \sum l_i=r}  \delta_1^{l_1,...,l_r}\, x_1^{l_1} \dots x_r^{l_r} &\equiv& x_1.(x_1+x_2)\dots (x_1+\dots+x_r)
\end{eqnarray}
We then let the formal derivation operator $\mi{log} \mg{\Pi}$ (rather than the algebra automorphism $\mg{\Pi}$ as in \S3.2) act on the test function $z$. {\it Without} resp. {\it with} ulterior derivation by $z$, we obtain these two expansions:
\begin{eqnarray}                         \label{vivi10}
\mg{\gop}_\ast(z)  \!\!&\!\!=\!\!&\!\!
\sum_{1\leq r}(-1)^{r\!-\!1}\!\!\! \sum_{ 0 \leq l_i < s_i}\!\!
 \,\delta^{l_1,...,l_r}\, \mathbf{Tan}^{s_1,..., s_r}(z)
\!\!  \prod_{1\leq i \leq r} \!\!
\frac{(s_i\!-\!1)!\,g_{\ast s_i\!-\!l_i\!+\!1}}{(s_i\!-\!l_i\!-\!1)!} 
  \\ [1.5 ex]                        \label{vivi11}
\mg{\gop}^{\prime}_\ast(z)  \!\!&\!\!=\!\!&\!\!
\sum_{1\leq r}(-1)^{r}\!\!\! \sum_{ 0 \leq l_i < s_i}\!\!
 \,\delta_1^{l_1,...,l_r}\, \mathbf{Tan}^{s_1,..., s_r}(z)
\!\!  \prod_{1\leq i \leq r}\!\!
\frac{(s_i\!-\!1)!\,g_{\ast s_i\!-\!l_i\!+\!1}}{(s_i\!-\!l_i\!-\!1)!}  
\end{eqnarray}
Of these expansions, the first is computationally more advantageous (it carries less summands) while the second is formally more appealing (its multitangents $\mi{Tan}^\bu$ have exactly the same total weight  $\sum s_j$ as the accompanying coefficient clusters). We may note that while it would be possible (but rather pointless) to produce similar expansions for all derivatives $\mg{\gop}_\ast^{(n)}$, nothing analogous exists for  the indefinite integrals 
$^{`\!}\mg{\gop}_\ast, ^{``\!}\mg{\gop}_\ast$\dots. In any case, and despite involving $g_\ast, \mg{\gop}_\ast$ rather
than  $g, \mg{\gop}$, the `symmetrical' expansions 
(\ref{vivi10}) and (\ref{vivi11}) 
 are theoretically simpler and more basic and computationally more efficient  than their `asymmetrical' counterparts 
 (\ref{vivi4}) and (\ref{vivi5}). 

%
\subsection{Parity separation.}

Formula (\ref{vivi10}) may be rewritten as
\begin{eqnarray*} 
\hspace{-3.ex }
\mg{\gop}_\ast(z)  \!\!&\!\!=\!\!&\!\!
\sum_{1\leq r}(-1)^{r\!-\!1}\!\!\! \sum^{\sum l_i=r-1}_{ 0 \leq l_i\, ,\, 1 \leq d_i}\!\!
 \delta^{l_1,...,l_r}\, \mathbf{Tan}^{d_1+l_1,..., d_r+l_r}\!(z)
\!\!  \prod_{1\leq i \leq r} \!\!
\frac{(d_i\!+\!l_i\!-\!1)!\,g_{\ast 1+d_i}}{(d_i\!-\!1)!}  
\end{eqnarray*}
Here, multitangents $\mi{Tan}^\bu$ of total ``degree"
$-r \p\sum(d_i\p l_i)=-1\p\sum d_i$ sit in front of coefficient clusters $\prod g_{\ast1+ d_i}$ of total ``degree"
$\sum d_i$. As {\it functions}, therefore, these multitangents always have the {\it same} parity as the accompanying coefficient clusters.


\subsection{From collectors to connectors.}
\noindent
{\bf From $\mg{\gop}$ to $\mg{\pi}=(\mg{\pi}_{\mr{no}}, \mg{\pi}_{\mr{so}} )$:}
\\
If we retain on the right-hand side of (\ref{vivi4}) the sole summands of weight\footnote{the `weight' in question is that of the coefficient clusters. But the weight of the accompanying multitangents (or, after reduction, of the multizeta-monotangent combinations) differs from the first only by one unit. } less than $s$ and call  
 $\mg{\gop}_{[s]}(z)$ the value of this finite sum, then, as $s$ goes to $\infty$ and for $K_\pm$ large enough,
the functions $\mg{\gop}_{[s]}(z)$ converge on half-planes $ \Im (z) > K_{+}$ resp.  $ \Im (z) <- K_{-}$  ,
to the northern resp. southern component of $\mg{\pi}$.

\noindent
{\bf From $\mg{\gop}_\ast$ to $\mg{\pi}_\ast=(\mg{\pi}_{\ast\mr{no}}, \mg{\pi}_{\ast\mr{so}} )$:}
\\
The `infinitesimal' $\mg{\gop}_\ast$ and $\mg{\pi}_\ast$ are formal, and generically divergent, series. But we may still consider the finite truncated sums $\mg{\gop}_{\ast [s]}(z)$ with their Fourier coefficients $\mg{\gop}_{\ast\omega [s]}$, and observe that, when $s\rightarrow$, these $\mg{\gop}_{\ast\omega [s]}$ do tend to finite limits $\mg{\gop}_{\ast\omega}$, which are the Fourier coefficients of the sought-after infinitesimal connector $\mg{\pi}_\ast$, or rather those of its northern and southern components.

\noindent
{\bf Remark:} 
Despite being very close to the {\it connectors}, the {\it collectors} differ from them in two fundamental respects: they are {\it not invariant} and they are {\it of one piece}. 

The {\it non-invariance} is fairly obvious for $\mg{\gop}$ taken in its {\it natural} multitangent expansion; but it remains true of 
$\mg{\gop}$ after its {\it natural} reduction to a monotangent expansion. Indeed, if we go from the `monotangential' 
$\mg{\gop}(z)$ to the power series $\mg{p}(z)$ by changing each $\mi{Te}^s(z)$
 to $z^{-s}$ and then formally Borel-transform $\mg{p}(z)$ to  $ \mg{\hat{p}}(\zeta)$ by changing $z^{-s}$ to 
 $\zeta^{s-1}/(s\!-\!1)!$, we end up with an entire function   $ \mg{\hat{p}}(\zeta)$ whose {\it only} invariant values
 are the ones it assumes on the set $\Omega:=2\pi i\,\doZ^\ast$. See \S4.1 {\it infra}.
 
 As for being {\it of one piece}, this is a property not so much of the collectors as of their constituent multitangents or monotangents\footnote{
 Things get more tangled in the {\it ramified case}, i.e. for diffeos with a contact order $p>1$. In that case, the multitangents still exist and still possess uniform determinations on each upper/lower half-plane of the $p$-ramified $z$-plane, but one can no longer go from one determination to the next by crossing the real axis between 
two consecutive integers.
 }, which are meromorphic over the whole of $\doC$, in complete contrast to the connectors, whose northern and southern components are usually completely unrelated: each one may a priori be {\it anything}.

%
\subsection{Reflexive and unitary diffeomorphisms.}
In \S3.4 we observed that in the expansion (\ref{vivi10}) of $\mg{\gop}_\ast$, 
coefficient clusters $\prod g_{\ast1+d_i}$ of even (resp. odd) total degree $\sum d_i$ 
accompany multitangents $\mi{Tan}^\bu$ that are even functions with real Fourier coefficients (resp. odd functions with purely imaginary Fourier coefficients). As a consequence, there is no simple condition on the coefficients $ g_{\ast1+d_i}$ of $g_\ast$ capable of ensuring that $\mg{\gop}_\ast $ be {\it odd}, whereas three elementary conditions may ensure that it be {\it even}, namely:
\\
(i) all coefficients $ g_{\ast1+d_i}$ of odd degree $d_i$ vanish and those of even degree are real
\\
(ii) all coefficients $ g_{\ast1+d_i}$ of even degree $d_i$ vanish and those of odd degree are purely imaginary
\\
(iii) all coefficients $ g_{\ast1+d_i}$ of even degree $d_i$ are real and those of odd degree are purely imaginary

No special significance attaches to case (ii), but the cases (i) and (iii) present interesting stability properties, with collectors and connectors inheriting the nature of $f$. This is an incentive for singling out the following three types of diffeos $f$ whose inverses  $f^{-1}$ either coincide with, or are analytically conjugate to, the  image of $f$ under an elementary involution:
\[\begin{array}{ccccccccccccccc}
\mi{reflexive} \!\!&:& \check{f}=f^{-1} &||& \mi{weakly\; reflexive} \!\!&:& \check{f} \stackrel{\mr{an.\,cj.}}{\sim} f^{-1} &
\\
\mi{unitary} \!\!&:& \bar{f}=f^{-1} &||& \mi{weakly\; unitary} \!\!&:& \bar{f} \stackrel{\mr{an.\,cj.}}{\sim} f^{-1} &
\\
\mi{counitary}\!\!&:& \check{\bar{f}}=f^{-1} &||& \mi{weakly\; counitary} \!\!&:& \check{\bar{f}} \stackrel{\mr{an.\,cj.}}{\sim} f^{-1} &
\end{array}\]
Here, $\bar{f}$ denotes the complex conjugate of $f$, and $\check{f}:=\sigma\circ f\circ \sigma$ with $\sigma(z)\equiv -z$.
Conjugation by $\tau$, with $\tau(z)\equiv i\,z$, clearly exchanges {\it unitary} and {\it counitary}, so that {\it weakly unitary} is equivalent to 
{\it weakly counitary}. Though {\it unitariness} seems a more natural notion, we shall work here with {\it counitariness}, which is better adapted to the correspondance $f\mapsto \mg{\pi}$ and enables us to take $f$ in standard form 
$f=l\circ g $.
\\     \noindent $\mathbf{P_1}$: 
 $f$ is reflexive iff the power series $f_\ast $ resp. $f^\ast$ are {\it even} resp. {\it odd}, in which case 
 $f_{\ast\pm}(-z)\equiv f_{\ast\mp}(z)$ and $ f^\ast_\pm(-z)\equiv -f^\ast_\mp(z)$.
 Likewise, $f$ is counitary iff the power series $f_\ast $ resp. $f^\ast$ are of the form  $f_{\ast\mr{re}}\circ \tau $ resp. 
 $\tau^{-1} \circ f_{\mr{re}}^\ast\circ \tau$ with
 real  $f_{\ast\mr{re}} $\,,\,$f^\ast_{\mr{re}}$, 
  in which case  $\bar{f}_{\ast\pm}(-z)\equiv f_{\ast\mp}(z)$ and $ \bar{f}^\ast_\pm(-z)\equiv -f^\ast_\mp(z)$.
\\     \noindent $\mathbf{P_2}$: If a standard $f$ is reflexive resp. counitary, then its conjugate 
$ l^{+\frac{1}{2}}\circ f \circ l^{-\frac{1}{2}} $ is of the standard form $f=l\circ g$ with reflexive resp. counintary factors $l$ and
$g:= l^{-\frac{1}{2}}\circ f \circ l^{-\frac{1}{2}} $.
\\     \noindent $\mathbf{P_3}$: If $f$ is (weakly or strictly) reflexive resp. counitary, then its connector $\mg{\pi}$ is (strictly) reflexive resp. counitary. This is geometrically obvious, from the relations $\mathbf{P_1}$ injected into the definition (\ref{ekam2}), but the remarkable fact is that the analytical procedure (\ref{vivi10}) also respects this conservation
of reflexivity or counitariness at every single step. Thus, if we apply it to the decomposition $f=l\circ g$ (as in $\mathbf{P_2}$) 
of a reflexive $f$, we have to do with an {\it even} infinitesimal generator $g_\ast$ that carries only coefficients $g_{\ast 1+d}$ of {\it even} degree $d$, and (\ref{vivi10}) automatically produces an {\it even} $\mg{p}_\ast$. The diffeo $g$ itself is of mixed parity, but its coefficients of $g_{\ast 1+d}$ of {\it odd} degree are fully determined by the earlier coefficients of {\it even} degree, and can thus be used in place of the  $g_{\ast 1+d}$. Either way, for reflexive diffeos the calculation of the invariants is a much more pleasant affair than for general diffeos, due to the drastic reduction in the mass of coefficients and (provided $f$ be real) to the realness of $\mg{\gop_\ast}$ and $\mg{\pi_\ast}$.
\\     \noindent $\mathbf{P_4}$: Conversely, any reflexive resp. counitary $\mg{\pi}$ is the invariant of some  
reflexive resp. counitary $f$. This follows from the {\it canonical synthesis} (see \S1.4) which, for $c$ real and large enough, automatically produces  diffeos $f_c$ of the required type.\footnote{ As pointed out to us by Reinhard Sch\"{a}ffke, this can also be deduced from the bifactorisation of $f$ in $\mathbf{P_5}$ below, provided we admit the existence of a pre-image $f$ for any given $\mg{\pi}$, which fact again follows from the canonical synthesis, but may also be established more directly.}
\\     \noindent $\mathbf{P_5}$: (Reinhard Sch\"{a}ffke). The product or quotient of two reflexive (res. unitary) diffeomorphisms is obviously conjugate to a reflexive (res. unitary) diffeomorphisms, but the converse is also true: any weakly reflexive (resp. unitary) $f$ can, for any consecutive integers $n_j$, be represented as a quotient of two strictly reflexive (resp. unitary) diffeos $f_j$:
\begin{equation*}
f:=f_1\circ f_2^{-1}  \quad \mi{with}\quad 
f(z):= z+1+o(1),\; 
f_j(z):= z+n_j+o(1), \;
n_1-n_2=1
\end{equation*}
and that too with explicit factors $f_j$:
\begin{eqnarray} \nonumber
f\;\; \mi{weakly\; reflexive}  
\quad\quad &||&\quad\quad
f \;\;\mi{weakly\; counitary}  \quad \quad
\\ [1.ex]   \label{weak1}
f_j:=      (^{\ast\!\!}f)   \circ l^{n_j}  \circ \sigma   \circ  (f^\ast)    \circ    \sigma
 &||&
f_j:= (^{\ast\!\!}f)   \circ l^{n_j}  \circ \sigma\circ  (\bar{f^\ast}) \circ \sigma 
\quad \quad
\\   \label{weak2}
\hspace{8.5 ex} =  f^{n_j} \circ (^{\ast\!\!}f)   \circ \sigma \circ  (f^\ast) \circ \sigma \!\!
&||&
\hspace{3.ex} 
=       f^{n_j}   \circ  (^{\ast\!\!}f) \circ \sigma \circ  (\bar{f^\ast}) \circ\sigma
\quad \quad
\\   \label{weak3}
\hspace{1. ex} = f^{n_j}  \circ h^{-1} \circ \sigma\circ h  \circ \sigma \quad
&||&
\hspace{3.ex}  = f^{n_j}  \circ h^{-1} \circ \sigma\circ \bar{h}  \circ \sigma \quad
\quad \quad
\end{eqnarray}
Indeed, the equivalent definitions (\ref{weak1}),  (\ref{weak2}), (\ref{weak3}) make it clear, respectively:
\\
\--- that $f_1, f_2$ are reflexive (resp. counitary); 
\\
\---  that $ f = f_1\circ f_2^{-1}$;
\\
\---  that $f_1, f_2$ are analytic.\footnote{The analytic $h$ in (\ref{weak3}) conjugates the weakly reflexive/counitary $f$ with a strictly reflexive/counitary $f_0$, i.e. $ h\circ f=f_0\circ h$. By definition, such a pair $h,f_0$ exists. 
We may note in passing that the factorisation $f=f_1\circ f_2^{-1} $ would still hold for complex  (in the reflexive case) or real (in the unitary case) values of $n_j$,
but in that case the above formulae break down ($f_1, f_2$ are no longer analytic) and we must take recourse to another, more involved construction.
}
\\     \noindent $\mathbf{P_6}$:
Piecing together all the above, we see that the commutative, non-associative\footnote{
$\mi{mix_c}(\mg{\pi_1},\mg{\pi_2})$ is doubly germinal: for a given $(\mg{\pi_1},\mg{\pi_2})$, it is defined for $c$ large enough, and for a given $c$ , it is defined for $(\mg{\pi_1},\mg{\pi_2})$ close enough to $(\mi{id},\mi{id})$.
} operation $\mi{mix_c}$:
\begin{equation} \label{mix}
\mr{mix_c}\;:\; (\mg{\pi_1},\mg{\pi_2}) \mapsto \mg{\pi}
:= \mg{\pi}_{f_{1,c}\circ f_{2,c}}
= \mg{\pi}_{f_{2,c}\circ f_{1,c}}
\end{equation}
(where $f_{j,c}$ stands for the $c$-canonical pre-image of $\mg{\pi}_j$) respects reflexivity and counitariness.

%

\section{Scalar invariants in terms of $f$.}
\subsection{The invariants $A_\omega$ as entire functions of $f$.}

Let $\mg{\pi}_\omega^{\pm}$ and  $\mg{\pi}_{\ast\omega}$ be the Fourier coefficient of the {\it connectors}, as defined in \S3.5
by weight-wise truncation of the {\it collectors} and passage to the limit:
\begin{eqnarray}          \label{catur1}
\hspace{-4.ex}
\!\mi{If}\, \p\Im(z)\!\gg\!1 
\!\!\!\!&:&\!\! \mg{\pi}^{\pm1}(z)=z+\sum_{\omega\in\Omega^{-}}  \mg{\pi}^{\pm}_\omega \,e^{-\omega z} 
\;; \;  \mg{\pi}_{\ast}(z)=\sum_{\omega\in\Omega^{-}}  \mg{\pi}_{\ast\omega} \,e^{-\omega z} 
\\ [0.7 ex]            \label{catur2}
\hspace{-4.ex}
\!\mi{If}\, \m\Im(z)\!\gg\!1 
\!\!\!\!&:&\!\! \mg{\pi}^{\pm1}(z)=z+\sum_{\omega\in\Omega^{+}}  \mg{\pi}^{\pm}_\omega \,e^{-\omega z} 
\;; \; \mg{\pi}_{\ast}(z)=\sum_{\omega\in\Omega^{+}}  \mg{\pi}_{\ast\omega} \,e^{-\omega z} 
\end{eqnarray}
In view of (\ref{vivi4})-(\ref{vivi5}) and (\ref{vivi10}), these Fourier coefficients are given by the convergent series\footnote{
Regarding the nature of their convergence, see the remark at the end of the present subsection.}
\begin{eqnarray}            \label{catur3}
\!\!\! \!\!\!\!\!\! \!\!\!
\mg{\pi}^{\pm}_\omega  &\!\!=\!\!&
\!\sum_{1\leq r}\!
\sum_{  {0 \leq l_i \atop 1\leq n_i }}^{n_i+l_i\leq s_i}
\!\!\!\!(-1)^{n\!-\!1}
 \, \delta^{\,l_1,...,\,l_r}_{n_1,...,n_r}\; \mathbf{Te}^{s_1,..., s_r}_\omega
\!\!  \prod_{1\leq i \leq r} \!\!
\frac{(s_i\!-\!1)!\,g^{\pm}_{n_i,s_i\!-\!l_i\!+\!1}}{(s_i\!-\!l_i\!-\!1)!} 
  \\ [1.5 ex]          \label{catur4}
\mg{\pi}_{\ast\omega}  &\!\!=\!\!&
\sum_{1\leq r}(-1)^{r\!-\!1}\!\! \sum_{ 0 \leq l_i < s_i}\!\!
 \,\delta^{l_1,...,l_r}\; \mathbf{Tan}^{s_1,..., s_r}_\omega
\!\!  \prod_{1\leq i \leq r} \!\!
\frac{(s_i\!-\!1)!\,g_{\ast s_i\!-\!l_i\!+\!1}}{(s_i\!-\!l_i\!-\!1)!} 
\end{eqnarray}
However, the need to define the alien operators $\Delta_\omega^\pm$ and $\Delta_\omega$ in uniform manner for all 
$\omega$ clashes with the need to associate within one and the same pair 
$(\mg{\pi}_{\mr{no}},\mg{\pi}_{\mr{so}}  )$ resp. $(\mg{\pi}^{-1}_{\mr{no}},\mg{\pi}^{-1}_{\mr{so}}  )$ 
northern and southern components originating from the same collector $\mg{\gop}$ or $\mg{\gop}^{-1}$.
This clash leads to a regrettable but unavoidable disharmony in the correspondance between the invariants 
$A_\omega^\pm$ and $A_\omega$, as defined from the resurgence equations, and the Fourier coefficients of the connectors, as derived from the collectors. This correspondance takes the form:
\[\begin{array}{llllllllllllllll}
\forall \omega\in\Omega^{-}
\;&: &\; A_\omega^{+}=\mg{\pi}_\omega^{+}
\;&;&\;  A_\omega^{-}= \mg{\pi}_\omega^{-}
\;&;&\;   +2\pi i\,A_\omega =\mg{\pi}_{\ast\omega}
\\ [1.5 ex]
\forall \omega\in\Omega^{+}
\;&: &\; A_\omega^{-} =\mg{\pi}_\omega^{+}
\;&;&\; A_\omega^{+}=  \mg{\pi}_\omega^{-}
\;&;&\;-2\pi i\,A_\omega =   \mg{\pi}_{\ast\omega} 
\end{array}\]
Alternatively, if we take the Borel transform $\mg{ \hat{p} }(\zeta)$ of the collector $\mg{p}(z)$ in its power series variant,
and observe that:
$$ \mg{\pi}^\pm_\omega\equiv \mi{sign(\Im \omega)}\,2\pi i\, \mg{\widehat{p}}^{\,\pm}(\omega) 
\quad \mi{and} \quad
 \mg{\pi}_{\ast\omega}\equiv \mi{sign(\Im \omega)}\,2\pi i\, \mg{\widehat{p}}_{\ast}(\omega) $$
then the correspondance assumes the following form, with yet another twist:
\[\begin{array}{llllllllllllllll}
\forall \omega\in\Omega^{-}
&: & A_\omega^{+} =-2\pi i\,\mg{\widehat{p}}^{\,+}(\omega)
&;&  A_\omega^{-} =-2\pi i\,\mg{\widehat{p}}^{\,-}(\omega)
&;&A_\omega =   -\mg{\widehat{p}_\ast}(\omega)
\\ [1.5 ex]
\forall \omega\in\Omega^{+}
&: & A_\omega^{-} =+2\pi i\,\mg{\widehat{p}}^{\,+}(\omega)
&;& A_\omega^{+} =  +2\pi i\,\mg{\widehat{p}}^{\,-}(\omega)
&;& A_\omega =  -\mg{\widehat{p}_\ast}(\omega)
\end{array}\]
{\bf Remark: nature of the convergence.}
\\
(i) Weight truncation: if we truncate (\ref{catur3}) to $\mg{\pi}^\pm_{\omega, [s]}$
or  (\ref{catur4}) to $\mg{\pi}_{\ast\omega, [s]}$ by retaining the sole terms of total weight $\leq s$ and then let $s$ go to $\infty$, we have convergence of the truncated series to the correct limits.
\\
(ii) If in (\ref{catur3}) we regroup all terms corresponding to identical (up to order) monomial 
clusters $\prod g_{s_i}^{k_i}$ into single blocks, then the series consisting of these blocks is always convergent, and may be regarded the Taylor expansion of  $\mg{\pi}^\pm_{\omega}$ viewed as an entire function of $g$.
\\
(iii) The same applies to  (\ref{catur4}) but {\it only} if we switch from the coefficients $g_{\ast s}$ (which have Gevrey-1 growth) to the $g_s$ (which have Gevrey-0 i.e. analytic growth) before forming the weight-homogeneous blocks. This does not contradict the point (i) above: when truncating (\ref{catur4}), there is no need of switching to the $g_s$, since truncation at weight $s$ is exactly the same whether performed relative to the coefficients of $g_\ast$ or those of $g$. 

%
\subsection{The case $\rho(f)\not=0$. Normalisation.}
For diffeos of the form $f(z)=z+1-\rho z^{-1}+\caO(z^{-2})$ with a non-vanishing `iterative residue' $\rho$, 
the defining relation (\ref{ekam1}) for the right and left iterators must be changed to
\begin{equation}    \label{ca1}
f_\pm^\ast(z)=\lim_{k\rightarrow \pm \infty} f^k(z)-k \pm\, \rho\,(c+\log |k|)
\end{equation}
with the normalisation constant $c$ as in \S2.6.
In the formal model, this leads to
\begin{equation}   \label{ca2}
\tilde{f}^\ast(z)=z+\rho\,(c+\log z)+ \caO(z^{-1})
\end{equation}
That apart, nothing changes and all the previous results and formulae still hold, including the explicit expansions 
(\ref{vivi4})-(\ref{vivi10}) and (\ref{catur3})-(\ref{catur4}),
provided we set $\mi{Ze}^1:=\gamma-c$ and normalise all multizetas and multitangents accordingly. As mentioned in \S2.6, the recommended choice is $c=\gamma$, since it amounts to setting $\mi{Ze}^1:=0$.


%
\subsection{The case $p(f)\not=1$. Ramification.}
Here again, the transition is straightforward. The `prepared' form (\ref{prep}) for the diffeo now carries fractional exponents
$s\in p^{-1}\,\doN^\ast$. As a consequence, the multiplicative $z$-plane and the convolutive $\zeta$-plane are now $p$-ramified, and so is the index set $\Omega$, which is embedded in the $\zeta$-plane. We still have one single collector
$\mg{\gop}$ (resp. $\mg{\gop}_\ast$), {\it ramified} yet {\it of one piece}, but $p$ distinct pairs of connectors,
$ \mg{\pi}=(\mg{\pi}_{\mr{no}}, \mg{\pi}_{\mr{so}})$ 
(resp. $ \mg{\pi}_\ast=(\mg{\pi}_{\ast\mr{no}}, \mg{\pi}_{\ast\mr{so}})$), separately {\it unramified} and mutually {\it unrelated}. The invariants $\mg{\pi}^{\pm}_\omega$ or  $\mg{\pi}_{\ast\omega}$ are still given by the familiar formulae
(\ref{catur3}), (\ref{catur4}), but with coefficients $\mi{Te}_\omega$ and  $\mi{Tan}_\omega$ that are best calculated by resurgent analysis, as in \S2.7, and are no longer finite sums of multizetas, even of ramified ones.

The transition to the most general case, with $(\rho,p)$ any element of $(\doC,\doN^\ast)$, follows exactly the same lines, and merely combines the partial adjustments of the present and preceding subsections.


%
\subsection{Growth properties of the invariants.}
{\bf Growth in $\omega$ for a given analytic $f$\,:}
\\
For a diffeo $f$ in prepared form (\ref{prep}), any majorisation of its coefficients easily translates into a majorisation of its invariants:
\begin{equation} \label{majmaj1}
\{\; |f_{[s]}| \leq c_0\, c_1^s  \;\} \Longrightarrow 
\{\; | A_\omega^\pm | \leq C_0\,C_1^{|\omega|} \; \}
\end{equation}
Rough estimates of $(C_0,C_1)$ in terms of  $(c_0,c_1)$ were given in [E2] and sharper ones in [B]. These results can be derived from a geometric analysis in the $z$-plane or from a resurgent analysis in the $\zeta$-plane. Things change, though, when we go over to the Gevrey case.

\noindent {\bf Growth in $\omega$ for a given $f$\,of Gevrey class:}
\\
Formal diffeos $f$ (in prepared form) of Gevrey class $\tau$ are easily 
shown to be stable under formal conjugations (also in prepared form) of the same Gevrey class. For $0< \tau$, the Gevrey class is non-analytic, and Gevrey conjugacy turns out to be strictly stronger than formal conjugacy if and only if
  $\tau<1$. 
This implies, for   $0<\tau<1$ ,  the existence of Gevrey conjugation invariants. These, however, can no longer be defined in the $z$-plane, since $f$ is purely formal and has no geometric realisation there. In the $\zeta$-plane, though, the Borel tranforms of the iterators 
$^{\ast\!\!}f$ and $f^\ast$ still exist (again, assuming $ \tau<1$); 
still extend to uniform analytic functions on $ \widetilde{\doC-2\pi i\,\doZ} $;
still verify the familiar resurgence equations 
(\ref{ekam22})-(\ref{ekam23}); and still unambigously defined define invariants $A_\omega^\pm$ and $A_\omega$, which are still given by the explicit expansions (\ref{catur3})-(\ref{catur4}). The only difference lies in the faster than exponential growth of $ \hat{f}^\ast(\zeta)$ and $^{\ast\!\!}\hat{f}(\zeta)$ as $|\zeta|\rightarrow \infty$, and in the faster than exponential growth of $A_\omega^{\pm}$ as $|\omega|\rightarrow\infty$. More precisly, for $0<\tau<1$, the 
earlier implication (\ref{majmaj1}) becomes\footnote{For details, see [E2], p 424.}:
\begin{equation} \label{majmaj2}
\{\; |f_{[s]}| \leq c_0\; c_1^s\; s^{\tau s}  \;\} \Longrightarrow 
\{\; | A_\omega^\pm | \leq C_0\;C_1^{|\omega|} \; \exp(C_2\,{|\omega|^\frac{1}{1-\tau}}) \}
\end{equation}

\noindent {\bf Growth in $f$ for a given  $\omega$\,:}
\\
We may now fix $\omega$ and ask how 
$A_\omega^{+}(f), A_\omega^{-}(f), A_\omega(f)$ behave as functions of $f$ or, to simplify, as entire functions
of any given coefficient $f_{[s]}$ $\,(s\geq 2$) relative to a prepared form (\ref{prep}). Unlike with the $\omega$-growth, 
there is little difference here between $A_\omega^{\pm}$ and $A_\omega$.
\\
(i) If $s>2$, all three entire functions $A_\omega^{+}(f_{[s]}), A_\omega^{-}(f_{[s]}), A_\omega(f_{[s]})$
have at most exponential growth in $ |f_{[s]}|^{\frac{1}{s-1}} $.
\\
(ii) If $s=2$, the corresponding coefficient  coincides up to sign with the iterative residue (i.e.
$f_{[2]}=-\rho$), and
the entire functions $A_\omega^{+}(\rho), A_\omega^{-}(\rho), A_\omega(\rho)$
have at most exponential growth in $ |\rho \log \rho| $. The result appears to be sharp.\footnote{See the argument
in \S8 of [BEE].}

These results are almost ``special cases" of the following statement: at any given point $\zeta_0$ on 
$\widetilde{\doC-\Omega}$, the Borel transform of the direct iterator assumes  a value $\hat{f}^\ast(\zeta_0)$
which, as an entire function of $f_{[s]}$, is exactly 
of exponential type in $ |f_{[s]}|^{\frac{1}{s-1}} $. This applies even for $s=2$. The difference between 
the cases $s\not=2$ and
  $s=2$ makes itself felt only when we move $\zeta_0$ to some point  $\omega_0$ located over $\Omega$, to investigate the leading singularity there and infer from it the value of the invariants. When $\rho=0$, the leading singularity
  in question  is a simple pole $a_{\omega_0} (\zeta-\omega_0)^{-1}$, but when $\rho\not=0$ it is of the form 
$a_{\omega_0} (\zeta-\omega_0)^{\rho\, \omega_0-1}/\Gamma( \rho\,\omega_0)$ and can be 
quite violent if $\rho$ has an imaginary part.

%
\subsection{Alternative computational strategies.}
\noindent
{\bf (i) Direct Fourier analysis in the multiplicative plane.}
\\
The methods amounts to calculating the limit:\footnote{If $\rho(f)\not=0$, the shift $l^{-k}$  should of course be replaced by
$l^{-k+(c+\log k)\rho}$, with $c=\gamma$ as recommended choice for the normalisation constant $c$. See \S2.6.
}
\begin{equation} \label{alter1}
A^{\mmpp\epsilon(\omega)}_\omega= \mg{\pi}^\pm_\omega 
=\lim_{k\rightarrow \pm \infty} \int_{z_0}^{1+z_0} 
\Big( l^{-k}\circ f^{2k}\circ l^{-k}(z)-z\Big)\,e^{\omega\,z} \,dz
\end{equation}
with $\epsilon(\omega):=\mi{sign}(\Im(\omega))$.
Although the parenthesised part of the integrand converges to $\mg{\pi}^\pm(z)-z $ for
$|\Im(z)|$ large enough, the above scheme, even after optimisation in the choice of $z_0$, 
is computationally costly (integral instead of series)
and inefficient (arithmetical convergence) as well as
 theoretically opaque (it sheds no light on the internal structure of the invariants
as functions of $f$).
\\

\noindent
{\bf (ii) Asymptotic coefficient analysis in the formal model. }
\\
The method starts with the inductive calculation of the first $N$ coefficients of the direct iterator $f^{\ast}(z)$ from its functional equation (\ref{iterat1}). One then switches to the Borel transform $\hat{f}^\ast(\zeta)$ and uses the method of {\it coefficient asymptotics}\footnote{
For a brief exposition of the method, see for ex. the section \S2.3 of {\it  
Power Series with sum-product Taylor coefficients and their resurgence algebra}, 
J. Ecalle and S. Sharma, Ed. Scuola Normale Superiore, Pisa, 2011. }
to derive the form of the two singularities\footnote{or of the $2\,p$ closest singularities when $p(f)\not=0$.
} closest to the origin (they are located over $\pm 2\pi i$). When applied to a parameter-free diffeo $f$ with proper optimising precautions, the method is superbly efficient for computing $A_{\pm 2\pi i}$, even for diffeos $f$ that are `large', i.e. distant from the identity. Thus, with $N$ taken in the region of 200 or 300, one typically gets $A_{\pm 2\pi i}$ with  100 exact digits or more, in less than half an hour of Maple time.
\\

\noindent
{\bf (iii) Resurgent analysis in the Poincar\'e plane. }
\\
The method is based on the resurgence equation (\ref{ekam23}) verified by the direct iterator $f^\ast$. But instead of interpreting that resurgence equation, as usual, in the highly ramified $\zeta$-plane, one performs a conformal transform
$\zeta\rightarrow \xi$ derived from  the classical modular function $\lambda$\,:
\begin{eqnarray} \label{poincare1}
\zeta\!&\!=\!&\!q(\xi):= -\log(1-\lambda(\xi))=-\log \lambda(-\frac{1}{\xi}) =16\sum_{n\,\mr{odd}} q_n \, e^{2\pi i \xi}
\quad \quad
\\   \label{poincare2}
q_n\!&\!:=\!&\!\sum_{d|n}\frac{1}{d}= \frac{1}{n}\sum_{d|n} d
\end{eqnarray}
 That comformal transform does three things:
\\
(*) it maps the Riemann surface $\widetilde{\doC -2\pi i\,\doZ}$ of the $\zeta$ variable 
uniformly 
onto the Poincar\'e half-plane $ \Im(\xi)>0$;
\\
(**) it changes the power series $\hat{f}^\ast(\zeta)$ with finite radius of convergence into a Fourier series
$\widehat{f}^\ast(\xi):=\hat{f}^\ast\circ q(\xi)$ that converges on the entire Poincar\'e half-plane.
\\
(***) it turns the alien operators into finite superpositions of post-composition operators \--- more precisely, post-composition  by simple homographies 
$h^\pm_{\omega,j}$ or $h^\pm_{\omega,j}$ with entire  coefficients:
\begin{eqnarray}   \label{ca3}
\Delta_\omega^\pm \widehat{\varphi}(\xi) &:=& 
\widehat{\varphi}\circ h^\pm_{\omega,1}(\xi)
-\widehat{\varphi}\circ h^\pm_{\omega,2}(\xi)
\\    \label{ca4}
\Delta_\omega \widehat{\varphi}(\xi) &:=& \sum_{1\leq j\leq 2^r}\, m_{\omega,j}\,
\widehat{\varphi}\circ h_{\omega,j}(\xi)
\hspace{ 5.ex} (r:=|\frac{\omega}{2\pi i}|, m_{\omega,j}\in \doQ)
\end{eqnarray}
The method is efficient enough for small values of $\omega$, but as $r:=|\frac{\omega}{2\pi i}|$ increases, the minima
\begin{eqnarray}   \label{ca5}
H^\pm(\omega) &:=&   \inf_{\Im(\xi)>0}  \{\, \Im(\xi)\, , \, \Im(h^\mp_{\omega,1}(\xi))\,,\, \Im(h^\mp_{\omega,2}(\xi))\}
\\   \label{ca6}
H(\omega) &:=&   \inf_{\Im(\xi)>0}  \{\, \Im(\xi)\, , \, \Im(h_{\omega,1}(\xi))\,, \dots,\, \Im(h_{\omega,2^r}(\xi))\}
\end{eqnarray}
rapidly decrease to zero, making it necessary to evaluate our Fourier series 
for $\widehat{f}^\ast(\xi)$ close to 
the boundary of their domain of convergence, i.e. the real axis, which of course is computationally costly. 
\\

\noindent
{\bf (iv) Explicit multizetaic expansions.}
\\
This method, to which the present paper was devoted, has the advantage of explicitness and theoretical transparency, expressing as it does the invariants in terms of universal transcendental constants (the multizetas) and of the diffeo's Taylor coefficients. It has the further advantage of handling large values of $\omega$ almost as efficiently as small ones. But the method's chief drawback would seem to be this: it involves expansions which converge very fast (faster than geometrically) once they reach `cruising speed', but which often take a damn long time to reach that speed. This is the case, not so much for $\omega$ large, but for $f$ large, i.e. for diffeos too distant from $\mi{id}$.


%
\subsection{Concluding remarks.}

\noindent
{\bf (i)  The invariants as autark functions. }
\\
Local, analytic, resonant vector fields $X$ ranging through a {\it fixed} formal conjugacy class, possess  holomorphic invariants $A_\omega$ which are {\it autark} functions of $X$, that is to say, of any given {\it free}\footnote{
i.e. of each coefficient that may freely vary without causing $X$ to leave its formal conjugacy class.} Taylor coefficient
of $X$. Autark functions, very informally, are entire functions whose asymptotic behaviour in every sector of exponential increase or decrease admits a complete description, with dominant exponential terms accompanied by  divergent-resurgent power series, which in turn verify a {\it closed} system of resurgence equations. Whether the invariants $A_\omega$ of diffeos are autark, too, seems likely but is yet unproved. Be that as it may, one would like to fully understand the asymptotic behaviour of $A_\omega$ as $f$ grows, or as any given coefficient or parameter in $f$ grows, since for
very  `large' diffeos $f$ the {\it direct} computation of the invariants would in any case be very costly.
\\

\noindent
{\bf (ii) Formal multizetas: dynamical vs arithmetical variants.}
\\
There exist several distinct but most probably equivalent notions of {\it arithmetical formal multizetas}, like the multizeta symbols subject to the two systems of so-called {\it quadratic multizeta relations}, or again to the {\it pentagonal, hexagonal} and {\it digonal relations}. But there also exists a demonstrably distinct and {\it weaker} notion of {\it dynamical formal multizetas} (and {\it multitangents}), by which we mean any system $\doS$ of 
scalar-valued multizeta symbols (resp. function-valued multitangent symbols) that, when inserted into the expansions (\ref{catur4}) (resp. (\ref{vivi4})) guarantees, first, the convergence of these expansions, and, second, the invariance of the $A_\omega$ (resp. $\mg{\pi}$) so produced. This immediately suggests a programme: to repeat for the dynamical multizetas what has been successfully done for their arithmetical counterparts, in particular to construct {\it explicit, complete and canonical systems of irreducibles}.
\\

{
\noindent
{\bf (iii) Abstract invariants. }
\\
Let $ \{\, ^{\doS\!\!}A_\omega,\, \omega\in \Omega \}$ be the system of `abstract' invariants induced by a system $\doS$ of dynamical multizetas as above. Since the system of natural invariants 
$ \{ A_\omega, \,\omega\in \Omega \}$ is complete, there necessarily exist conversion formulae of the form:
\begin{equation}   \label{ca10}
^{\doS\!\!}A_{\omega_0}=\sum_{1\leq r}\sum_{\omega_1+...\omega_r=\omega_0}
H_\doS^{\omega_1,...,\omega_r} A_{\omega_1}\dots A_{\omega_r}
\end{equation}
that respect the basic $\omega$-gradation and carry interesting `universal' structure constants $H_\doS^\bu$.
These constants ought to be of particular significance in the case of the system $\doS_0$ of `rational' dynamical multizetas which is analogous, on the dynamical side, to the canonical system of `rational'\footnote{they become rational, of course, only after an homogeneous rescaling that amounts to setting $\pi:=1$. } multizetas on the arithmetical side.

\subsection{Some historical background.}

\noindent
{\bf (i) Identity-tangent diffeos in holomorphic dynamics. }
\\
The iteration of one-dimensional analytic mappings -- whether local or global; identity-tangent or not -- has a long history going back a century or more.  Fatou, for one, knew about the analytic classes of identity-tangent diffeos and had formed a clear, geometry-based idea of their invariants. The subject then when into something of a hibernation, until the advent of high-power computation, which brought about an explosive revival of holomorphic dynamics, one- and many-dimensional. For the specific subject of analytic invariants, however, the main impetus for renewal came from an unexpected quarter: resurgent analysis.
\\

\noindent
{\bf (ii) Identity-tangent diffeos and resurgent analysis. }
\\
The fact is that identity-tangent diffeos possess generically divergent but always resurgent iterators and fractional iterates, with an interesting, non-linear pattern of resurgence or self-reproduction at the singular points in the Borel plane, and it was in the process of sorting  out these phenomena that resurgence theory was born, and later applied to general local objects and much else. In a sense, this involved a retreat from dynamics proper, since it meant focusing on the Borel plane, where the key dynamic notions of trajectory, fixed point etc admit no simple interpretation. For the invariants $A_\omega$, however, the shift in focus brought a definite advantage, since in the Borel plane these invariants are automatically {\it localised} and 
{\it isolated} (they appear as coefficients of the leading singularities over the point $\omega$) whereas in the multiplicative plane they are {\it diffuse} and {\it intertwined} (they make themselves felt only collectively and indirectly, via Stokes phenomena and the like, and the only way to isolate them is by Fourier analysis of type (\ref{alter1}), which is but a  half-hearted way of doing what Borel analysis does neatly and efficiently). This applies not just to identity-tangent diffeos, but to a huge range of local objects and equations. It also works in both directions: in that of ``analysis", i.e. calculating and investigating the invariants of a given object; and in that of ``synthesis", i.e. starting from an admissible system of prospective `invariants' and constructing an object of which they are the actual invariants. And it has to be said that in both directions resurgence theory performs rather better than geometry. It leads in particular to a privileged or ``canonical" synthesis, a notion which eludes geometry.
\\

\noindent
{\bf (iii) Identity-tangent diffeos and the resuscitation of multizetas. }
Multizetas, especially of length 2, were first considered by Euler as an isolated curiosity, and later fell into a protracted oblivion for want of applications. They resurfaced only in the late 1970s and early 1980s, precisely in the context of holomorphic dynamics and identity-tangent diffeos, as {\it the} transcendental ingredient in the make-up of their invariants. Ten years later, the multizetas started cropping up in half a dozen, largely unconnected contexts: braid groups and knot theory; Feynman diagrams;  Galois theory; mixed Tate motives; arithmetical dimorphy; ARI/GARI and the flexion structure, etc.
At the moment, all these strands are in the process of merging or at least cross-fertilising, and constitute a vibrant, active field of research.
\\

\noindent
{\bf (iv) Identity-tangent diffeos and the actual computation of their invariants. }
\\
The sections of  [E2] devoted to the invariants of identity-tangent diffeos were written with no computational applications in mind, and no attempt was made to optimise the bounds in formulae such as  (\ref{majmaj1}) or (\ref{majmaj2}). On the contrary, the PhD thesis [B], which revisits the subject 30 years on, lays its main emphasis on these neglected aspects and provides effective Maple programmes for the computations of the invariants; it also offers copious asides on the algebraic aspects of multitangents, which largely, but not exactly, mirror those of multizetas.




%


\section{Tables.}
\subsection{Multitangents: from symmetrel to alternal.}
We express $\mi{Tan}^\bu$ in terms of $\mi{Te}^\bu$ according to the defining relations
$$
\{\mr{Tan}^\bu = (\mr{logmu}.\mr{Te}^\bu) \circ (E^\bu-\mg{1}^\bu)\}
\Longleftrightarrow
\{ (E^\bu-\mg{1}^\bu)\circ \mr{Tan}^\bu =
\mr{Te}^\bu \circ (E^\bu-\mg{1}^\bu)\}
$$
\[\begin{array}{lllll}
\mr{Tan}^{n_1}  \!\!&\! \!&\! =  \mr{Te}^{n_1}  
\\ [1.5 ex]
\mr{Tan}^{n_1,n_2}  \!\!&\! \!&\! =
\frac{1}{2}\mr{Te}^{n_1, n_2}
-\frac{1}{2}\mr{Te}^{n_2, n_1}
\\ [1.5 ex]
\mr{Tan}^{n_1,n_2,n_3}  \!\!&\!\!&\!  =
\frac{1}{3}\mr{Te}^{n_1,n_2,n_3}
-\frac{1}{6}\mr{Te}^{n_1,n_3,n_2}
-\frac{1}{6}\mr{Te}^{n_2,n_1,n_3}
-\frac{1}{6}\mr{Te}^{n_2,n_3,n_1}
-\frac{1}{6}\mr{Te}^{n_3,n_1,n_2}
\\ [1. ex]  \!&\!\!&\!  
+\frac{1}{3}\mr{Te}^{n_3,n_2,n_1}
-\frac{1}{6}\mr{Te}^{n_1+n_3,n_2}
+\frac{1}{12}\mr{Te}^{n_1,n_2+n_3}
+\frac{1}{12}\mr{Te}^{n_1+n_2,n_3}
\\ [1.ex]  \!&\!\!&\!  
+\frac{1}{12}\mr{Te}^{n_3,n_1+n_2}
+\frac{1}{12}\mr{Te}^{n_2+n_3,n_1}
-\frac{1}{6}\mr{Te}^{n_2,n_1+n_3}
\\ [1.5 ex]
\mr{Tan}^{n_1,n_2,n_3,n_4}  \!\!\!\!&\!\!&\! = 
\frac{1}{4}\mr{Te}^{n_1, n_2, n_3, n_4}
-\frac{1}{12}\mr{Te}^{n_1, n_2, n_4, n_3}
-\frac{1}{12}\mr{Te}^{n_1, n_3, n_2, n_4}
-\frac{1}{12}\mr{Te}^{n_1, n_3, n_4, n_2}
\\ [1. ex]
\!\!\!\!&\! \!&\!  
-\frac{1}{12}\mr{Te}^{n_1, n_4, n_2, n_3}
+\frac{1}{12}\mr{Te}^{n_1, n_4, n_3, n_2}
-\frac{1}{12}\mr{Te}^{n_2, n_1, n_3, n_4}
+\frac{1}{12}\mr{Te}^{n_2, n_1, n_4, n_3}
\\  [1. ex]
\!\!\!\!&\! \!&\!  
-\frac{1}{12}\mr{Te}^{n_2, n_3, n_1, n_4}
-\frac{1}{12}\mr{Te}^{n_2, n_3, n_4, n_1}
+\frac{1}{12}\mr{Te}^{n_2, n_4, n_1, n_3}
+\frac{1}{12}\mr{Te}^{n_2, n_4, n_3, n_1}
\\  [1. ex]
\!\!\!\!&\! \!&\!  
-\frac{1}{12}\mr{Te}^{n_3, n_1, n_2, n_4}
-\frac{1}{12}\mr{Te}^{n_3, n_1, n_4, n_2}
+\frac{1}{12}\mr{Te}^{n_3, n_2, n_1, n_4}
+\frac{1}{12}\mr{Te}^{n_3, n_2, n_4, n_1}
\\  [1. ex]
\!\!\!\!&\! \!&\!  
-\frac{1}{12}\mr{Te}^{n_3, n_4, n_1, n_2}
+\frac{1}{12}\mr{Te}^{n_3, n_4, n_2, n_1}
-\frac{1}{12}\mr{Te}^{n_4, n_1, n_2, n_3}
+\frac{1}{12}\mr{Te}^{n_4, n_1, n_3, n_2}
\\  [1. ex]
\!\!\!\!&\! \!&\!  
+\frac{1}{12}\mr{Te}^{n_4, n_2, n_1, n_3}
+\frac{1}{12}\mr{Te}^{n_4, n_2, n_3, n_1}
+\frac{1}{12}\mr{Te}^{n_4, n_3, n_1, n_2}
-\frac{1}{4}\mr{Te}^{n_4, n_3, n_2, n_1}
\\  [1. ex]
\!\!\!\!\!\!\!\!&\! \!&\!   \!\!\!\!\!\!\!\!
+\frac{1}{12}\mr{Te}^{n_1, n_2, n_3+n_4}
-\frac{1}{12}\mr{Te}^{n_1, n_3, n_2+n_4}
-\frac{1}{12}\mr{Te}^{n_2, n_3, n_1+n_4}
+\frac{1}{12}\mr{Te}^{n_2, n_4, n_1+n_3}
\\  [1. ex]
\!\!\!\!\!\!\!\!&\! \!&\!   \!\!\!\!\!\!\!\!
-\frac{1}{12}\mr{Te}^{n_3, n_1, n_2+n_4}
+\frac{1}{12}\mr{Te}^{n_3, n_2,n_1+n_4}
+\frac{1}{12}\mr{Te}^{n_4, n_2, n_1+n_3}
-\frac{1}{12}\mr{Te}^{n_4, n_3, n_1+n_2}
\\  [1. ex]
\!\!\!\!\!\!\!\!&\! \!&\!   \!\!\!\!\!\!\!\!
+\frac{1}{12}\mr{Te}^{n_1, n_2+n_3, n_4}
-\frac{1}{12}\mr{Te}^{n_1, n_2+n_4, n_3}
-\frac{1}{12}\mr{Te}^{n_2, n_1+n_3, n_4}
+\frac{1}{12}\mr{Te}^{n_2, n_1+n_4, n_3}
\\  [1. ex]
\!\!\!\!\!\!\!\!&\! \!&\!   \!\!\!\!\!\!\!\! 
-\frac{1}{12}\mr{Te}^{n_3, n_1+n_4, n_2}
+\frac{1}{12}\mr{Te}^{n_3, n_2+n_4, n_1}
+\frac{1}{12}\mr{Te}^{n_4, n_1+n_3, n_2}
-\frac{1}{12}\mr{Te}^{n_4, n_2+n_3, n_1}
\\  [1. ex]
\!\!\!\!\!\!\!\!&\! \!&\!   \!\!\!\!\!\!\!\!
+\frac{1}{12}\mr{Te}^{n_1+n_2, n_3, n_4}
-\frac{1}{12}\mr{Te}^{n_1+n_3, n_2, n_4}
-\frac{1}{12}\mr{Te}^{n_1+n_3, n_4, n_2}
+\frac{1}{12}\mr{Te}^{n_2+n_4, n_1, n_3}
\\  [1. ex]
\!\!\!\!\!\!\!\!&\! \!&\!   \!\!\!\!\!\!\!\!
+\frac{1}{12}\mr{Te}^{n_2+n_4, n_3, n_1}
-\frac{1}{12}\mr{Te}^{n_1+n_4, n_2, n_3}
+\frac{1}{12}\mr{Te}^{n_1+n_4, n_3, n_2}
-\frac{1}{12}\mr{Te}^{n_3+n_4, n_2, n_1}
\\  [1. ex]
\!\!\!\!\!\!\!\!&\! \!&\!   \!\!\!\!\!\!\!\!
+\!\frac{1}{24}\mr{Te}^{n_1+n_2, n_3+n_4}
\!-\! \frac{1}{12}\mr{Te}^{n_1+n_3, n_2+n_4}
\!+\! \frac{1}{12}\mr{Te}^{n_2+n_4, n_1+n_3}
\!-\! \frac{1}{24}\mr{Te}^{n_3+n_4, n_1+n_2}
\\ [1.5 ex]
\mr{Tan}^{n_1,\dots,n_5}  \!\!\!\!&\!\!&\! = \;\; \mi{ 540\;\; summands}.
\\ [1.5 ex]
\mr{Tan}^{n_1,\dots,n_6}  \!\!\!\!&\!\!&\! = \;\; \mi{ 3688\;\; summands}.
\\ [1.5 ex]
\mr{Tan}^{n_1,\dots,n_7}  \!\!\!\!&\!\!&\! = \;\; \mi{ 47292\;\; summands}.
\end{array}\]
\subsection{Parity/imparity of alternal multitangents.}
We begin by comparing the number of summands in the monotangent reductions 
$\mi{red_1(Te^\bu)}$ and $\mi{red_1(Tan^\bu)}$ 
(resp. $\mi{red_2(Te^\bu)}$ and $\mi{red_2(Tan^\bu)}$)
of $\mi{Te^\bu}$ and $\mi{Tan^\bu}$ 
{\it before} (resp. {\it after}) symmetrel linearisation of the resulting multizetas. N.B. A further reduction
$\mi{red_3(Te^\bu)}$ and $\mi{red_3(Tan^\bu)}$, corresponding to a complete decomposition of the multizeta into {\it arithmetical irreducibles}, would yield even fewer summands.
\[\begin{array}{rrrrrrrrrrrrrrrrrr}
(n_1,...,n_r) &||&\# \mr{red_1(Te^\bu)}\!\!\! &|&\# \mr{red_1(Tan^\bu)} \!\!\!   &||&\# \mr{red_2(Te^\bu)} \!\!\! &|&\# \mr{red_2(Tan^\bu)} \!\!\! \!\!
\\[1.7 ex]
 (5,7,4)        &||&  34        &|&   40      &||&   30      &|& 17
\\[1. ex]
 (5,7,14)      &||& 133         &|&  124       &||&  106       &|&  65
\\[1. ex]
 (5, 7, 4, 5)   &||&   148       &|&  141       &||&   177      &|&  80
 \\[1. ex]
 (8, 11, 4, 9)   &||&  580        &|&  679       &||&  1127       &|&  454 
\\[1. ex]
 (8, 11, 7, 12)  &||&   824       &|&   741      &||&  1154       &|&  452
\\[1. ex]
 (4,5,4,5,4)     &||&   42       &|&   54      &||&   389      &|&  98
\\[1. ex]
 (3,4,5,6,7)     &||&   455       &|&  874       &||&  2748       &|&  760
\end{array}\]
The following six examples of multitangent reduction (of type $\mi{red}_2$) illustrate the phenomenon of {\it parity separation} in 
$\mi{Tan}^\bu$, as opposed to $\mi{Te}^\bu$.
\\

\noindent {\bf Example 1\,:} $\mr{Te}^{2,6,4}$ is neither even nor odd.
\[\begin{array}{rllll}
\mr{Te}^{2,6,4}(z) \!\!\!&\!=\!&\!\!\!   \sum_{2 \leq m\leq 6} \mr{Teze}^{2,6,4}_m\, \mr{Te}^{m}(z)\;\mi{with}
\\ [1.9 ex]
\mr{Teze}^{2,6,4}_1\!\!\!&\!=\!&\!\!\! 
+20\, \mr{Ze}^{11}
+56\, \mr{Ze}^{9, 2}
-70\, \mr{Ze}^{8, 3}
-112\, \mr{Ze}^{3, 8}
+54\, \mr{Ze}^{7, 4}
+42\, \mr{Ze}^{4, 7}
\\  \!\!\!&&\!\!\!
-20\, \mr{Ze}^{6, 5}
-20\, \mr{Ze}^{5, 6}
\\    \!\!\!&\!=\!&\!\!\! 
\, 0
\\
\mr{Teze}^{2,6,4}_2 \!\!\!&\!=\!&\!\!\! 
+14\,\mr{Ze}^{10}
+35\,\mr{Ze}^{8, 2}
+56\,\mr{Ze}^{2, 8}
-40\,\mr{Ze}^{7, 3}
-28\,\mr{Ze}^{3, 7}
+35\,\mr{Ze}^{6, 4}
\\  \!\!\!&&\!\!\!
+39\,\mr{Ze}^{4, 6}
-32\,\mr{Ze}^{5, 5}
\\
\mr{Teze}^{2,6,4}_3 \!\!\!&\!=\!&\!\!\! 
+8\,\mr{Ze}^{9}
+20\,\mr{Ze}^{7, 2}
+14\,\mr{Ze}^{2, 7}
-20\,\mr{Ze}^{6, 3}
-22\,\mr{Ze}^{3, 6}
+8\,\mr{Ze}^{5, 4}
+8\,\mr{Ze}^{4, 5}
\\
\mr{Teze}^{2,6,4}_4 \!\!\!&\!=\!&\!\!\! 
+5\,\mr{Ze}^{8}
+10\,\mr{Ze}^{6, 2}
+11\,\mr{Ze}^{2, 6}
-8\,\mr{Ze}^{5, 3}
-8\,\mr{Ze}^{3, 5}
+6\,\mr{Ze}^{4, 4}
\\
\mr{Teze}^{2,6,4}_5 \!\!\!&\!=\!&\!\!\! 
+2\,\mr{Ze}^{7}
+4\,\mr{Ze}^{5, 2}
+4\,\mr{Ze}^{2, 5}
-2\,\mr{Ze}^{4, 3}
-2\,\mr{Ze}^{3, 4}
\\
\mr{Teze}^{2,6,4}_6 \!\!\!&\!=\!&\!\!\! 
+\mr{Ze}^{6}
+\mr{Ze}^{4, 2}
+\mr{Ze}^{2, 4}
\end{array}\]

\noindent {\bf Example 2\,:} $\mi{Tan}^{2,6,4}$ is  even since $\mi{deg}(\mi{Tan}^{2,6,4} )\!=\!2\p6\p4\m3=9=\mi{odd} $.
\[\begin{array}{rllll}
\mr{Tan}^{2,6,4}(z) \!\!\!&\!=\!&\!\!\!   
\mr{Tanze}^{2,6,4}_2\, \mr{Te}^{2}(z)+\mr{Tanze}^{2,6,4}_4\, \mr{Te}^{4}(z)+\mr{Tanze}^{2,6,4}_6\, \mr{Te}^{6}(z)\;\mi{with}
\\ [1.5 ex]
\mr{Tanze}^{2,6,4}_2 \!\!\!&\!=\!&\!\!\! 
+5\,\mr{Ze}^{10}
-\frac{7}{3}\,\mr{Ze}^{8, 2}
+\frac{56}{3}\,\mr{Ze}^{2, 8}
-40\,\mr{Ze}^{7, 3}
-28\,\mr{Ze}^{3, 7}
+9\,\mr{Ze}^{6, 4}
\\  [0.5 ex]    \!\!\!&&\!\!\! 
+13\,\mr{Ze}^{4, 6}
-32\,\mr{Ze}^{5, 5}
\\  [0.5 ex]
\mr{Tanze}^{2,6,4}_4 \!\!\!&\!=\!&\!\!\! 
+3\,\mr{Ze}^{8}
+\frac{8}{3}\,\mr{Ze}^{6, 2}
+\frac{11}{3}\,\mr{Ze}^{2, 6}
-8\,\mr{Ze}^{5, 3}
-8\,\mr{Ze}^{3, 5}
+2\,\mr{Ze}^{4, 4}
\\   [0.5 ex]
\mr{Tanze}^{2,6,4}_6 \!\!\!&\!=\!&\!\!\! 
+\frac{2}{3}\,\mr{Ze}^{6}
+\frac{1}{3}\,\mr{Ze}^{4, 2}
+\frac{1}{3}\,\mr{Ze}^{2, 4}
\end{array}\]

\noindent {\bf Example 3\,:} $\mr{Te}^{2,7,4}(z)$ is neither even nor odd.
\[\begin{array}{rllll}
\mr{Te}^{2,7,4}(z) \!\!\!&\!=\!&\!\!\! \sum_{2 \leq m\leq 7} \mr{Teze}^{2,7,4}_m\, \mr{Te}^{m}(z)\;\mi{with}
\\ [1.5 ex]
\mr{Teze}^{2,7,4}_1\!\!\!&\!=\!&\!\!\! 
+30\,\mr{Ze}^{12}
+84\,\mr{Ze}^{10, 2}
-112\,\mr{Ze}^{9, 3}
-168\,\mr{Ze}^{3, 9}
+112\,\mr{Ze}^{8, 4}
\\    [0.5 ex]  \!\!\!&&\!\!\! 
+84\,\mr{Ze}^{4, 8}
-104\,\mr{Ze}^{7, 5}
-112\,\mr{Ze}^{5, 7}
+100\,\mr{Ze}^{6, 6}
\\    [0.5 ex]  \!\!\!&\!=\!&\!\!\! 
\, 0
\\  [0.5 ex]
\mr{Teze}^{2,7,4}_2 \!\!\!&\!=\!&\!\!\! 
+20\,\mr{Ze}^{11}
+56\,\mr{Ze}^{9, 2}
+84\,\mr{Ze}^{2, 9}
-70\,\mr{Ze}^{8, 3}
-56\,\mr{Ze}^{3, 8}
+54\,\mr{Ze}^{7,4}
\\    [0.5 ex]  \!\!\!&&\!\!\! 
+56\,\mr{Ze}^{4, 7}
-20\,\mr{Ze}^{6, 5}
-20\,\mr{Ze}^{5, 6}
\\   [0.5 ex]
\mr{Teze}^{2,7,4}_3 \!\!\!&\!=\!&\!\!\! 
+14\,\mr{Ze}^{10}
+35\,\mr{Ze}^{8, 2}
+28\,\mr{Ze}^{2, 8}
-40\,\mr{Ze}^{7, 3}
-42\,\mr{Ze}^{3, 7}
+35\,\mr{Ze}^{6, 4}
\\     [0.5 ex]  \!\!\!&&\!\!\! 
+35\,\mr{Ze}^{4, 6}
-32\,\mr{Ze}^{5, 5}
\\    [0.5 ex]
\mr{Teze}^{2,7,4}_4 \!\!\!&\!=\!&\!\!\! 
+8\,\mr{Ze}^{9}
+20\,\mr{Ze}^{7, 2}
+21\,\mr{Ze}^{2, 7}
-20\,\mr{Ze}^{6, 3}
-20\,\mr{Ze}^{3, 6}
+8\,\mr{Ze}^{5, 4}
+8\,\mr{Ze}^{4, 5}
\\    [0.5 ex]
\mr{Teze}^{2,7,4}_5 \!\!\!&\!=\!&\!\!\! 
+5\,\mr{Ze}^{8}
+10\,\mr{Ze}^{6, 2}
+10\,\mr{Ze}^{2, 6}
-8\,\mr{Ze}^{5, 3}
-8\,\mr{Ze}^{3, 5}
+6\,\mr{Ze}^{4, 4}
\\  [0.5 ex]
\mr{Teze}^{2,7,4}_6 \!\!\!&\!=\!&\!\!\! 
+2\,\mr{Ze}^{7}
+4\,\mr{Ze}^{5, 2}
+4\,\mr{Ze}^{2, 5}
-2\,\mr{Ze}^{4, 3}
-2\,\mr{Ze}^{3, 4}
\\    [0.5 ex]
\mr{Teze}^{2,7,4}_7 \!\!\!&\!=\!&\!\!\! 
+\mr{Ze}^{6}
+\mr{Ze}^{4, 2}
+\mr{Ze}^{2, 4}
\end{array}\]

\noindent {\bf Example 4\,:} $\mi{Tan}^{2,7,4}$ is  even since $\mi{deg}(\mi{Tan}^{2,7,4} )\!=\!2\p7\p4\m3\!=\!10=\mi{odd} $.
\[\begin{array}{rllll}
\mr{Tan}^{2,7,4}(z) \!\!\!&\!=\!&\!\!\!   
 \mr{Tanze}^{2,7,4}_3\, \mr{Te}^{3}(z)
  +\mr{Tanze}^{2,7,4}_5\, \mr{Te}^{5}(z)
   +\mr{Tanze}^{2,7,4}_7\, \mr{Te}^{7}(z)\;\mi{with}
\\ [1.5 ex]
\mr{Tanze}^{2,7,4}_1  \!\!\!&\!=\!&\!\!\! 
+36\,\mr{Ze}^{12}
+84\,\mr{Ze}^{10, 2}
-112\,\mr{Ze}^{9, 3}
-168\,\mr{Ze}^{3, 9}
+56\,\mr{Ze}^{8, 4}
\\   [0.5 ex]  \!\!\!&&\!\!\! 
+28\,\mr{Ze}^{4, 8}
-104\,\mr{Ze}^{7, 5}
-112\,\mr{Ze}^{5, 7}
+\frac{100}{3}\,\mr{Ze}^{6, 6}
\\    \!\!\!&\!=\!&\!\!\! 
\, 0
\\
\mr{Tanze}^{2,7,4}_3    \!\!\!&\!=\!&\!\!\! 
+11\,\mr{Ze}^{10}
+\frac{49}{3}\,\mr{Ze}^{8, 2}
+\frac{28}{3}\,\mr{Ze}^{2, 8}
-40\,\mr{Ze}^{7, 3}
-42\,\mr{Ze}^{3, 7}
+\frac{35}{3}\,\mr{Ze}^{6, 4}
\\   [0.5 ex]  \!\!\!&&\!\!\! 
+\frac{35}{3}\,\mr{Ze}^{4, 6}
-32\,\mr{Ze}^{5, 5}
\\   [0.5 ex]
\mr{Tanze}^{2,7,4}_5    \!\!\!&\!=\!&\!\!\! 
+\frac{10}{3}\,\mr{Ze}^{8}
+\frac{10}{3}\,\mr{Ze}^{6, 2}
+\frac{10}{3}\,\mr{Ze}^{2, 6}
-8\,\mr{Ze}^{5, 3}
-8\,\mr{Ze}^{3, 5}
+2\,\mr{Ze}^{4, 4}
\\  [0.5 ex]
\mr{Tanze}^{2,7,4}_7  \!\!\!&\!=\!&\!\!\! 
+\frac{2}{3}\,\mr{Ze}^{6}
+\frac{1}{3}\,\mr{Ze}^{4, 2}
+\frac{1}{3}\,\mr{Ze}^{2, 4}
\end{array}\]

\noindent {\bf Example 5\,:} $\mi{Tan}^{2,5,2,4}$ is  even since $\mi{deg}(\mi{Tan}^{2,5,2,4})\!=\!2\p5\p2\p4\m4\!=\!9=\mi{odd} $.
\[\begin{array}{rllll}
\mr{Tan}^{2,5,2,4}(z) \!\!\!&\!=\!&\!\!\!   \mr{Tanze}^{2,5,2,4}_2\, \mr{Te}^{2}(z)+\mr{Tanze}^{2,5,2,4}_4\, \mr{Te}^{4}(z)\;\mi{with}
\\ [1.5 ex]
\mr{Tanze}^{2,5,2,4}_2     \!\!\!&\!=\!&\!\!\! 
+\frac{14}{3}\mr{Ze}^{9, 2}
+\frac{2}{3}\mr{Ze}^{2, 9}
+8\,\mr{Ze}^{8, 3}
+\frac{4}{3}\mr{Ze}^{3, 8}
+\frac{50}{3}\mr{Ze}^{7, 4}
-\mr{Ze}^{4, 7}
+5\,\mr{Ze}^{6, 5}
\\   [0.5 ex]   \!\!\!&&\!\!\! 
+\frac{20}{3}\mr{Ze}^{5, 6}
+\frac{40}{3}\mr{Ze}^{7,2, 2}
+\frac{65}{3}\mr{Ze}^{2, 7, 2}
-\frac{20}{3}\mr{Ze}^{6, 3, 2}
-\frac{20}{3}\mr{Ze}^{6, 2, 3}
-\frac{10}{3}\mr{Ze}^{3, 6, 2}
\\    [0.5 ex]  \!\!\!&&\!\!\! 
+\frac{70}{3}\mr{Ze}^{2, 6, 3}
+\frac{20}{3}\mr{Ze}^{3, 2, 6}
+\frac{20}{3}\mr{Ze}^{2, 3, 6}
+\frac{20}{3}\mr{Ze}^{5, 4, 2}
+\frac{32}{3}\mr{Ze}^{5, 2, 4}
+7\,\mr{Ze}^{4, 5, 2}
\\   [0.5 ex]
   \!\!\!&&\!\!\! 
+\frac{32}{3}\mr{Ze}^{2, 5, 4}
+5\,\mr{Ze}^{2, 4, 5}
+6\,\mr{Ze}^{4, 4, 3}
+4\,\mr{Ze}^{4, 3, 4}
+2\,\mr{Ze}^{3, 4, 4}
\\  [0.5 ex]
   \!\!\!&&\!\!\! 
-32\,\mr{Ze}^{5, 3, 3}
-12\,\mr{Ze}^{3, 5, 3}
\\  [0.5 ex]
\mr{Tanze}^{2,5,2,4}_4      \!\!\!&\!=\!&\!\!\! 
+\mr{Ze}^{7, 2}
-\frac{1}{3}\mr{Ze}^{2,7}
+2\,\mr{Ze}^{6, 3}
+\frac{2}{3}\mr{Ze}^{3, 6}
+\frac{8}{3}\mr{Ze}^{5, 4}
+\frac{8}{3}\mr{Ze}^{5, 2, 2}
+\frac{7}{3}\mr{Ze}^{2, 5, 2}
\\  [0.5 ex]
   \!\!\!&&\!\!\! 
-\frac{2}{3}\mr{Ze}^{3, 4, 2}
+2\,\mr{Ze}^{2, 4, 3}
+\frac{4}{3}\mr{Ze}^{3, 2, 4}
+\frac{4}{3}\mr{Ze}^{2, 3, 4}
\end{array}\]

\noindent {\bf Example 6\,:} $\mi{Tan}^{2,3,2,5}$ is  even since $\mi{deg}(\mi{Tan}^{2,3,2,5})\!=\!2\p3\p2\p5\m4\!=\!8\!=\mi{even} $.
\[\begin{array}{rllll}
\mr{Tan}^{2,3,2,5}(z) \!\!\!&\!=\!&\!\!\!   
  \mr{Tanze}^{2,3,2,5}_3\, \mr{Te}^{3}(z)
+\mr{Tanze}^{2,3,2,5}_5\, \mr{Te}^{5}(z)\;\mi{with}
\\ [1.5 ex]
\mr{Tanze}^{2,3,2,5}_1     \!\!\!&\!=\!&\!\!\! 
+\frac{2}{3}\mr{Ze}^{9, 2}
+\frac{10}{3}\mr{Ze}^{2, 9}
-11\,\mr{Ze}^{8, 3}
-\frac{8}{3}\mr{Ze}^{3, 8}
-7\,\mr{Ze}^{7, 4}
+2\,\mr{Ze}^{4, 7}
-\frac{10}{3}\mr{Ze}^{5, 6}
\\  [0.5 ex]
  \!\!\!&&\!\!\! 
-10\,\mr{Ze}^{7, 2, 2}
+10\,\mr{Ze}^{2, 7, 2}
-\frac{25}{3}\mr{Ze}^{6, 3, 2}
-\frac{40}{3}\mr{Ze}^{6, 2, 3}
+10\,\mr{Ze}^{3, 6, 2}
-\frac{25}{3}\mr{Ze}^{2, 6, 3}
\\  [0.5 ex]
  \!\!\!&&\!\!\! 
-20\,\mr{Ze}^{3, 2, 6}
-20\,\mr{Ze}^{2, 3, 6}
-2\,\mr{Ze}^{5, 4, 2}
-8\,\mr{Ze}^{5, 2, 4}
+5\,\mr{Ze}^{4, 5, 2}
+10\,\mr{Ze}^{2, 5, 4}
\\  [0.5 ex]
  \!\!\!&&\!\!\! 
-10\,\mr{Ze}^{4, 2, 5}
+5\,\mr{Ze}^{2, 4, 5}
+20\,\mr{Ze}^{5, 3, 3}
+30\,\mr{Ze}^{3, 5, 3}
+10\,\mr{Ze}^{3, 3, 5}
-6\,\mr{Ze}^{4, 4, 3}
\\  [0.5 ex]
  \!\!\!&&\!\!\! 
-15\,\mr{Ze}^{4, 3, 4}
-6\,\mr{Ze}^{3, 4, 4}
\\  [0.5 ex]
\mr{Tanze}^{2,3,2,5}_3     \!\!\!&\!=\!&\!\!\! 
+\frac{2}{3}\mr{Ze}^{7, 2}
+\frac{4}{3}\mr{Ze}^{2, 7}
-\frac{2}{3}\mr{Ze}^{6, 3}
-\frac{2}{3}\mr{Ze}^{3, 6}
+\frac{1}{3}\mr{Ze}^{5, 4}
+\mr{Ze}^{4, 5}
-\frac{2}{3}\mr{Ze}^{5, 2, 2}
\\   [0.5 ex]
  \!\!\!&&\!\!\! 
+\frac{17}{3}\mr{Ze}^{2, 5, 2}
-\mr{Ze}^{4, 3, 2}
-4\,\mr{Ze}^{4, 2, 3}
+2\,\mr{Ze}^{3, 4, 2}
+2\,\mr{Ze}^{2, 4, 3}
-4\,\mr{Ze}^{3, 2, 4}
\\  [0.5 ex]
  \!\!\!&&\!\!\! 
-\mr{Ze}^{2, 3, 4}
+4\,\mr{Ze}^{3, 3, 3}
\\  [0.5 ex]
\mr{Tanze}^{2,3,2,5}_5    \!\!\!&\!=\!&\!\!\! 
+\frac{1}{3}\mr{Ze}^{5, 2}
+\frac{1}{3}\mr{Ze}^{2, 5}
+\mr{Ze}^{2, 3, 2}
\end{array}\]

\subsection{The invariants as entire functions of $f$: the general case.}
We write down, up to weight 10 inclusively, the expansion of the collector $\mg{\gop_\ast}$ in terms of
the $g_\ast$. We assume $p(f)=1$ but impose no restriction on $\rho(f)\equiv -g_{\ast2}$. In this and all further examples,
we order the terms according to their total weight and, within a given total weight, we start with the lowest monotangents.
\\

\noindent
{\bf Example 1}:\, $\mg{\gop_\ast}$ up to weight 10 for $f=l\circ g$ with $g_\ast(z)=\sum_{1\leq d} g_{\ast 1+d}z^{-d}$.
\begin{eqnarray*}
+\mathbf{Te}^{1}
\!\!\!&\!\!\!&\!\!\!
\big[ \mr{g}_{\ast2 } \big]
+\mathbf{Te}^{2}\,\big[\mr{g}_{\ast3 }\big]
+\mathbf{Te}^{3}\,\big[\mr{g}_{\ast4 }\big]
+\mathbf{Te}^{4}\,\big[\mr{g}_{\ast5 }\big]
+\mathbf{Te}^{2}\,\big[
6\,\zeta(3)\,\mr{g}_{\ast2 }\,\mr{g}_{\ast4 }
-6\,\zeta(3)\,\mr{g}_{\ast3 }^2
\big]
\\
+\mathbf{Te}^{5}
\!\!\!&\!\!&\!\!\!
\big[\mr{g}_{\ast6 }\big]
+\mathbf{Te}^{3}\,\big[
6\,\zeta(3)\,\mr{g}_{\ast2 }\,\mr{g}_{\ast5 }
-6\,\zeta(3)\,\mr{g}_{\ast3 }\,\mr{g}_{\ast4 }
\big[
+\mathbf{Te}^{6}\,\big[\mr{g}_{\ast7 }\big]
\\
+\mathbf{Te}^{2}
\!\!\!&\!\!&\!\!\!
\big[
30\,\zeta(5)\,\mr{g}_{\ast4 }^2
-\frac{5}{2}\,\zeta(5)\,\mr{g}_{\ast2 }^4
+10\,\zeta(5)\,\mr{g}_{\ast2 }\,\mr{g}_{\ast6 }
-40\,\zeta(5)\,\mr{g}_{\ast3 }\,\mr{g}_{\ast5 }
\big]
\\
+\mathbf{Te}^{3}
\!\!\!&\!\!&\!\!\!
\big[
\frac{4}{3}\,\zeta(2)^2\,\mr{g}_{\ast2 }\,\mr{g}_{\ast3 }^2
-\frac{4}{3}\,\zeta(2)^2\,\mr{g}_{\ast2 }^2\,\mr{g}_{\ast4 }
\big]
+\mathbf{Te}^{4}\big[
3\,\zeta(3)\,\mr{g}_{\ast4 }^2
+\frac{1}{4}\,\zeta(3)\,\mr{g}_{\ast2 }^4
-10\,\zeta(3)\,\mr{g}_{\ast3 }\,\mr{g}_{\ast5 }
\\
\!\!\!&\!\!&\!\!\!
+7\,\zeta(3)\,\mr{g}_{\ast2 }\,\mr{g}_{\ast6 }
\big]
+\mathbf{Te}^{5}
\big[
-\frac{2}{3}\,\zeta(2)\,\mr{g}_{\ast2 }\,\mr{g}_{\ast3 }^2
+\frac{2}{3}\,\zeta(2)\,\mr{g}_{\ast2 }^2\,\mr{g}_{\ast4 }
\big]
+\mathbf{Te}^{7}\,\big[\mr{g}_{\ast8 }\big]
\\
+\mathbf{Te}^{2}
\!\!\!&\!\!&\!\!\!
\big[
36\,\zeta(3)^2\,\mr{g}_{\ast3 }^3
-\frac{32}{5}\,\zeta(2)^3\,\mr{g}_{\ast3 }^3
+18\,\zeta(3)^2\,\mr{g}_{\ast5 }\,\mr{g}_{\ast2 }^2
+\frac{48}{5}\,\zeta(2)^3\,\mr{g}_{\ast2 }\,\mr{g}_{\ast3 }\,\mr{g}_{\ast4 }
\\
\!\!\!&\!\!&\!\!\!
 -54\,\zeta(3)^2\,\mr{g}_{\ast2 }\,\mr{g}_{\ast3 }\,\mr{g}_{\ast4 }
-\frac{16}{5}\,\zeta(2)^3\,\mr{g}_{\ast5 }\,\mr{g}_{\ast2 }^2
\big]
+\mathbf{Te}^{3}\,\big[
20\,\zeta(5)\,\mr{g}_{\ast4 }\,\mr{g}_{\ast5 }
+10\,\zeta(5)\,\mr{g}_{\ast2 }\,\mr{g}_{\ast7 }
\\
\!\!\!&\!\!&\!\!\!
-30\,\zeta(5)\,\mr{g}_{\ast3 }\,\mr{g}_{\ast6 }
-5\,\zeta(5)\,\mr{g}_{\ast2 }^3\,\mr{g}_{\ast3 }
\big]
+\mathbf{Te}^{4}\,\big[
-\frac{1}{5}\,\zeta{2}^2\,\mr{g}_{\ast3 }^3
-\frac{21}{10}\,\zeta(2)^2\,\mr{g}_{\ast2 }^2\,\mr{g}_{\ast5 }
\\
\!\!\!&\!\!&\!\!\!
+\frac{23}{10}\,\zeta(2)^2\,\mr{g}_{\ast2 }\,\mr{g}_{\ast3 }\,\mr{g}_{\ast4 }
\big]
+\mathbf{Te}^{5}\,\big[
8\,\zeta(3)\,\mr{g}_{\ast2 }\,\mr{g}_{\ast7 }
-12\,\zeta(3)\,\mr{g}_{\ast3 }\,\mr{g}_{\ast6 }
+4\,\zeta(3)\,\mr{g}_{\ast4 }\,\mr{g}_{\ast5 }
\\
\!\!\!&\!\!&\!\!\!
+\zeta(3)\,\mr{g}_{\ast2 }^3\,\mr{g}_{\ast3 }
\big]
+\mathbf{Te}^{6}\,\big[
-\frac{1}{3}\,\zeta(2)\,\mr{g}_{\ast3 }^3
+\frac{3}{2}\,\zeta(2)\,\mr{g}_{\ast2 }^2\,\mr{g}_{\ast5 }
-\frac{7}{6}\,\zeta(2)\,\mr{g}_{\ast2 }\,\mr{g}_{\ast3 }\,\mr{g}_{\ast4 }
\big]
\\
+\mathbf{Te}^{8}
\!\!\!&\!\!&\!\!\!
\big[\mr{g}_{\ast9}\big]
+\mathbf{Te}^{2}\,\big[
210\,\zeta(7)\,\mr{g}_{\ast4 }\,\mr{g}_{\ast6 }
-140\,\zeta(7)\,\mr{g}_{\ast5 }^2
-84\,\zeta(7)\,\mr{g}_{\ast3 }\,\mr{g}_{\ast7 }
+14\,\zeta(7)\,\mr{g}_{\ast2 }\,\mr{g}_{\ast8 }
\\
\!\!\!&\!\!&\!\!\!
-\frac{133}{3}\,\zeta(7)\,\mr{g}_{\ast2 }^3\,\mr{g}_{\ast4 }
+\frac{133}{3}\,\zeta(7)\,\mr{g}_{\ast2 }^2\,\mr{g}_{\ast3 }^2
\big]
+\mathbf{Te}^{3}\,\big[
36\,\zeta(3)^2\,\mr{g}_{\ast3 }^2\,\mr{g}_{\ast4 }
-9\,\zeta(3)^2\,\mr{g}_{\ast2 }\,\mr{g}_{\ast4 }^2
\\
\!\!\!&\!\!&\!\!\!
+21\,\zeta(3)^2\,\mr{g}_{\ast2 }^2\,\mr{g}_{\ast6 }
+\frac{3}{4}\,\zeta(3)^2\,\mr{g}_{\ast2 }^5
-\frac{32}{5}\,\zeta(2)^3\,\mr{g}_{\ast3 }^2\,\mr{g}_{\ast4 }
-\frac{64}{15}\,\zeta(2)^3\,\mr{g}_{\ast2 }^2\,\mr{g}_{\ast6 }
\\
\!\!\!&\!\!&\!\!\!
+\frac{32}{3}\,\zeta(2)^3\,\mr{g}_{\ast2 }\,\mr{g}_{\ast3 }\,\mr{g}_{\ast5 }
-48\,\zeta(3)^2\,\mr{g}_{\ast2 }\,\mr{g}_{\ast3 }\,\mr{g}_{\ast5 }
\big]
+\mathbf{Te}^{4}\,\big[
45\,\zeta(5)\,\mr{g}_{\ast4 }\,\mr{g}_{\ast6 }
-20\,\zeta(5)\,\mr{g}_{\ast5 }^2
\\
\!\!\!&\!\!&\!\!\!
-36\,\zeta(5)\,\mr{g}_{\ast3 }\,\mr{g}_{\ast7 }
+11\,\zeta(5)\,\mr{g}_{\ast2 }\,\mr{g}_{\ast8 }
-\frac{10}{3}\,\zeta(5)\,\mr{g}_{\ast2 }^3\,\mr{g}_{\ast4 }
-\frac{25}{6}\,\zeta(5)\,\mr{g}_{\ast2 }^2\,\mr{g}_{\ast3 }^2
\big]
\\
+\mathbf{Te}^{5}
\!\!\!&\!\!&\!\!\!
\big[
\frac{10}{3}\,\zeta(2)^2\,\mr{g}_{\ast2 }\,\mr{g}_{\ast3 }\,\mr{g}_{\ast5 }
-\frac{2}{5}\,\zeta(2)^2\,\mr{g}_{\ast3 }^2\,\mr{g}_{\ast4 }
-\frac{44}{15}\,\zeta(2)^2\,\mr{g}_{\ast2 }^2\,\mr{g}_{\ast6 }
\big]
+\mathbf{Te}^{6}\,\big[
9\,\zeta(3)\,\mr{g}_{\ast2 }\,\mr{g}_{\ast8 }
\\
\!\!\!&\!\!&\!\!\!
-14\,\zeta(3)\,\mr{g}_{\ast3 }\,\mr{g}_{\ast7 }
+5\,\zeta(3)\,\mr{g}_{\ast4 }\,\mr{g}_{\ast6 }
+\frac{1}{2}\,\zeta(3)\,\mr{g}_{\ast2 }^2\,\mr{g}_{\ast3 }^2
+2\,\zeta(3)\,\mr{g}_{\ast2 }^3\,\mr{g}_{\ast4 }
\big]
\\
+\mathbf{Te}^{7}
\!\!\!&\!\!&\!\!\!
\big[
\frac{8}{3}\,\zeta(2)\,\mr{g}_{\ast2 }^2\,\mr{g}_{\ast6 }
-\frac{5}{3}\,\zeta(2)\,\mr{g}_{\ast2 }\,\mr{g}_{\ast3 }\,\mr{g}_{\ast5 }
-\,\zeta(2)\,\mr{g}_{\ast3 }^2\,\mr{g}_{\ast4 }
\big]
+\mathbf{Te}^{9}\,\big[\mr{g}_{\ast10}\big]
\end{eqnarray*}

\subsection{The invariants as entire functions of $f$: the reflexive case.}
As in Example 1, we write down the expansion of the collector $\mg{\gop_\ast}$ in terms of
the $g_\ast$, but this time for a reflexive $g$.  We still assume $p(f)=1$ and reflexivity automatically implies
 $\rho(f)\equiv -g_{\ast2}\equiv 0$. There being fewer coefficients $g_{\ast s}$, we reach weight 13.
\\

\noindent
{\bf Example 2}:\, $\mg{\gop_\ast}$ up to weight 13 for $f=l\circ g$ with $g_\ast(z)=\sum_{1\leq d} g_{\ast 1+2d}z^{-2d}$.
\begin{eqnarray*}
+\mathbf{Te}^{2}
\!\!\!&\!\!&\!\!\!
\big[ \mr{g}_{\ast3} \big]
+\mathbf{Te}^{4}
 \big[ \mr{g}_{\ast5 }\big]
+\mathbf{Te}^{2}\;
\big[ \m6\,\zeta(3)\,\mr{g}_{\ast3}^2 \big]
+\mathbf{Te}^{6}\;
\big[\mr{g}_{\ast7 }\big]
+\mathbf{Te}^{2}\;
\big[\m40\,\zeta(5)\,\mr{g}_{\ast3 }\,\mr{g}_{\ast5 }\big]
\\
+\mathbf{Te}^{4}
\!\!\!&\!\!\!&\!\!\!
\big[\!-\!10\,\zeta(3)\,\mr{g}_{\ast3 }\,\mr{g}_{\ast5 }\big]
\!+\mathbf{Te}^{2}\;
\big[\,36\,\zeta(3)^2\,\mr{g}_{\ast3 }^3
\!-\!\frac{32}{5}\,\zeta(2)^3\,\mr{g}_{\ast3 }^3  \big]
\!+\mathbf{Te}^{4}\;
\big[\!-\!\frac{1}{5}\,\zeta(2)^2\,\mr{g}_{\ast3 }^3\big]
\\
+\mathbf{Te}^{6}
\!\!\!&\!\!\!&\!\!\!
\big[ \!-\!\frac{1}{3}\,\zeta(2)\,\mr{g}_{\ast3 }^3 \big]
+\mathbf{Te}^{8}\;\big[ \mr{g}_{\ast9}\big]
+\mathbf{Te}^{2}\;
\big[\!-\!140\,\zeta(7)\,\mr{g}_{\ast5 }^2
            \!-\!84\,\zeta(7)\,\mr{g}_{\ast3 }\,\mr{g}_{\ast7 }\big]
\\
+\mathbf{Te}^{4}
\!\!\!&\!\!\!&\!\!\!
     \big[\m20\,\zeta(5)\,\mr{g}_{\ast5 }^2
             \m36\,\zeta(5)\,\mr{g}_{\ast3 }\,\mr{g}_{\ast7 }\big]
+\mathbf{Te}^{6}\;
\big[\m14\,\zeta(3)\,\mr{g}_{\ast3 }\,\mr{g}_{\ast7 }]
\\
+\mathbf{Te}^{2} 
\!\!\!&\!\!\!&\!\!\!
\big[\m\frac{15648}{175}\,\zeta(2)^4\,\mr{g}_{\ast3 }^2\,\mr{g}_{\ast5 }
            -80\,\zeta(6,2)\,\mr{g}_{\ast3 }^2\,\mr{g}_{\ast5 }
            +800\,\zeta(3)\,\zeta(5)\,\mr{g}_{\ast3 }^2\,\mr{g}_{\ast5 }\big]
\\
+\mathbf{Te}^{4}
\!\!\!&\!\!\!&\!\!\!
            \big[ \m \frac{272}{21}\,\zeta(2)^3\,\mr{g}_{\ast3 }^2\,\mr{g}_{\ast5 }
             +80\,\zeta(3)^2\,\mr{g}_{\ast3 }^2\,\mr{g}_{\ast5 }\big]
+\mathbf{Te}^{6}\;
\big[\m \frac{34}{15}\,\zeta(2)^2\,\mr{g}_{\ast3 }^2\,\mr{g}_{\ast5 }\big]
\\
+\mathbf{Te}^{8}
\!\!\!&\!\!\!&\!\!\!
\big[\m \frac{5}{3}\,\zeta(2)\,\mr{g}_{\ast3 }^2\,\mr{g}_{\ast5 }\big]
+\mathbf{Te}^{10}\;\big[\mr{g}_{\ast11 }\big]
+\mathbf{Te}^{2}\;
\big[\m144\,\zeta(9)\,\mr{g}_{\ast3 }\,\mr{g}_{\ast9 }
             \m1008\,\zeta(9)\,\mr{g}_{\ast5 }\,\mr{g}_{\ast7 }
             \\   
\!\!\!&\!\!\!&\!\!\!
               \m 210\,\zeta(9)\,\mr{g}_{\ast3 }^4
               \m216\,\zeta(3)^3\,\mr{g}_{\ast3 }^4
               \p\frac{576}{5}\,\zeta(3)\zeta(2)^3\,\mr{g}_{\ast3 }^4\big]
\\
+\mathbf{Te}^{4}
\!\!\!&\!\!\!&\!\!\!
\big[\frac{18}{5}\,\zeta(3)\,\zeta(2)^2\,\mr{g}_{\ast3 }^4
                            \p14\,\zeta(7)\,\mr{g}_{\ast3 }^4
                            \m78\,\zeta(7)\,\mr{g}_{\ast3 }\,\mr{g}_{\ast9 }
                           \m210\,\zeta(7)\,\mr{g}_{\ast5 }\,\mr{g}_{\ast7 }\big]
\\
+\mathbf{Te}^{6}
\!\!\!&\!\!\!&\!\!\!
\big[\,6\,\zeta(2)\,\zeta(3)\,\mr{g}_{\ast3 }^4
                \m\frac{10}{3}\,\zeta(5)\,\mr{g}_{\ast3 }^4
                 \m28\,\zeta(5)\,\mr{g}_{\ast5 }\,\mr{g}_{\ast7 }
                 \m44\,\zeta(5)\,\mr{g}_{\ast3 }\,\mr{g}_{\ast9 }\big]
\\
+\mathbf{Te}^{8}
\!\!\!&\!\!\!&\!\!\!
\big[\m18\,\zeta(3)\,\mr{g}_{\ast3 }\,\mr{g}_{\ast9 }\big]
+\mr{Te}^{2}\;
\big[\m168\,\zeta(8,2)\,\mr{g}_{\ast3 }^2\,\mr{g}_{\ast7 }
             -280\,\zeta(8,2)\,\mr{g}_{\ast3 }\,\mr{g}_{\ast5 }^2
\\  
\!\!\!&\!\!\!&\!\!\!
            \m\frac{125056}{385}\,\zeta(2)^5\,\mr{g}_{\ast3 }\,\mr{g}_{\ast5 }^2
             \m\frac{375168}{1925}\,\zeta(2)^5\,\mr{g}_{\ast3 }^2\,\mr{g}_{\ast7 }
             \p1760\,\zeta(5)^2\,\mr{g}_{\ast3 }\,\mr{g}_{\ast5 }^2
\\  
\!\!\!&\!\!\!&\!\!\!
             \p1056\,\zeta(5)^2\,\mr{g}_{\ast3 }^2\,\mr{g}_{\ast7 }
             \p 3360\,\zeta(3)\,\zeta(7)\,\mr{g}_{\ast3 }\,\mr{g}_{\ast5 }^2
             \p 2016\,\zeta(3)\,\zeta(7)\,\mr{g}_{\ast3 }^2\,\mr{g}_{\ast7 }\big]
\\
+\mathbf{Te}^{4}
\!\!\!&\!\!\!&\!\!\!
\big[1080\,\zeta(3)\,\zeta(5)\,\mr{g}_{\ast3 }^2\,\mr{g}_{\ast7 }
             \p \frac{23824}{175}\,\zeta(2)^4\,\mr{g}_{\ast3 }^2\,\mr{g}_{\ast7 }
             \p180\,\zeta(6,2)*\mr{g}_{\ast3}^2\,\mr{g}_{\ast7 }
 \\
\!\!\!&\!\!\!&\!\!\!
              \p\frac{6544}{525}\,\zeta(2)^4\,\mr{g}_{\ast3 }\,\mr{g}_{\ast5 }^2
               \p100\,\zeta(6,2)\,\mr{g}_{\ast3 }\,\mr{g}_{\ast5 }^2
               \p200\,\zeta(3)\,\zeta(5)\,\mr{g}_{\ast3 }\,\mr{g}_{\ast5 }^2\big]
\\
+\mathbf{Te}^{6}
\!\!\!&\!\!\!&\!\!\!
\big[\m\frac{3064}{105}\,\zeta(2)^3\,\mr{g}_{\ast3 }^2\,\mr{g}_{\ast7 }
                      \m140\,\zeta(3)^2\,\mr{g}_{\ast3 }^2\,\mr{g}_{\ast7 }
                  \p\frac{88}{21}\,\zeta(2)^3\,\mr{g}_{\ast3 }\,\mr{g}_{\ast5 }^2\big]
\\
+\mathbf{Te}^{8}
\!\!\!&\!\!\!&\!\!\!
\big[\frac{8}{15}\,\zeta(2)^2\,\mr{g}_{\ast3 }\,\mr{g}_{\ast5 }^2
             \m\frac{39}{5}\,\zeta(2)^2\,\mr{g}_{\ast3 }^2\,\mr{g}_{\ast7 }\big]
 \\            
+\mathbf{Te}^{10}
\!\!\!\!\!&\!\!\!&\!\!
\big[\m\frac{2}{3}\,\zeta(2)\,\mr{g}_{\ast3 }\,\mr{g}_{\ast5 }^2
                        \m4\,\zeta(2)\,\mr{g}_{\ast3 }^2\,\mr{g}_{\ast7 })
+\mathbf{Te}^{12}\; \big[\mr{g}_{\ast13 }\big]
\end{eqnarray*}
\subsection{The invariants as entire functions of $f$: one-parameter cases.}
{\bf Example 3}:\, $\mg{\gop_\ast}$ up to weight 12 for $f=l\circ g$ with $g_(z)=z+g_{2}\,z^{-1}$.
\begin{eqnarray*}
+\mathbf{Te}^{1}
   \!\!\!&\!\!\!&\!\!\!
\,\mr{g}_{2}
+\mathbf{Te}^{3}\,\big[\frac{1}{2}\big]\,\mr{g}_{2}^2
+\mathbf{Te}^{2}\,\big[3\,\zeta(3)\big]\,\mr{g}_{2}^3
+\mathbf{Te}^{5}\,\big[\frac{1}{2}\big]\,\mr{g}_{2}^3
+\mathbf{Te}^{2}\,\big[10\,\zeta(5)\big]\,\mr{g}_{2}^4
\\
+\mathbf{Te}^{3}
   \!\!\!&\!\!\!&\!\!\!
\big[\m\frac{2}{3}\,\zeta(2)^2\big]\,\mr{g}_{2}^4
+\mathbf{Te}^{4}
\big[\frac{9}{2}\,\zeta(3)\big]\,\mr{g}_{2}^4
\!+\mathbf{Te}^{5}\,\big[\frac{1}{3}\,\zeta(2)\big]\,\mr{g}_{2}^4
\!+\mathbf{Te}^{7}\,\big[\frac{7}{12}\big]\,\mr{g}_{2}^4
\!+\mathbf{Te}^{2}
\big[\frac{77}{2}\,\zeta(7)\big]\,\mr{g}_{2}^5
\\
+\mathbf{Te}^{3}
   \!\!\!&\!\!\!&\!\!\!
\big[9\,\zeta(3)^2\m\frac{32}{15}\,\zeta(2)^3\big]\mr{g}_{2}^5
+\mathbf{Te}^{4}
\big[16\,\zeta{5}\big]\,\mr{g}_{2}^5
+\mathbf{Te}^{5}\,\big[\m\frac{22}{15}\,\zeta(2)^2\big]\,\mr{g}_{2}^5
+\mathbf{Te}^{6}\,\big[\frac{15}{2}\,\zeta(3)\big]\,\mr{g}_{2}^5
\\
+\mathbf{Te}^{7}
   \!\!\!&\!\!\!&\!\!\!
\big[\frac{4}{3}\,\zeta(2)\big]\,\mr{g}_{2}^5
+\mathbf{Te}^{9}\,\big[\frac{2}{3}\big]\,\mr{g}_{2}^5
+\mathbf{Te}^{2}
\big[151\,\zeta(9)\big]\,\mr{g}_{2}^6
+\mathbf{Te}^{3}\,\big[54\,\zeta(3)\,\zeta(5)]-\frac{14758}{2625}\,\zeta(2)^4
\\
   \!\!\!&\!\!\!&\!\!\!
+3\,\zeta(6,2)\big]\,\mr{g}_{2}^6
+\mathbf{Te}^{4}\,\big[\frac{271}{4}\,\zeta(7)-6\,\zeta(3)\,\zeta(2)^2\big]\,\mr{g}_{2}^6
+\mathbf{Te}^{5}
\big[27\,\zeta(3)^2-\frac{1052}{175}\,\zeta(2)^3\big]\,\mr{g}_{2}^6
\\
+\mathbf{Te}^{6}
   \!\!\!&\!\!\!&\!\!\!
\big[\frac{55}{2}\,\zeta(5)+5\,\zeta(2)\,\zeta(3)\big]\,\mr{g}_{2}^6
+\mathbf{Te}^{7}\,\big[-\frac{134}{75}\,\zeta(2)^2\big]\,\mr{g}_{2}^6
+\mathbf{Te}^{8}
\big[\frac{193}{16}\,\zeta(3)\big]\,\mr{g}_{2}^6
\\
+\mathbf{Te}^{9}
   \!\!\!&\!\!\!&\!\!\!
\big[\frac{53}{15}\,\zeta(2)\big]\,\mr{g}_{2}^6
+\mathbf{Te}^{11}\,\big[\frac{13}{20}\big]\,\mr{g}_{2}^6
\end{eqnarray*}
\noindent 
{\bf Example 4}:\, $\mg{\gop_\ast}$ up to weight 12 for $f=l\circ g$ with 
$ g(z)=z\,\Big[1+2\,g_{\ast2} z^{-2} \Big]^{\frac{1}{2}}$.
\begin{eqnarray*}
+\mathbf{Te}^{1}
   \!\!\!&\!\!&\!\!\!
\big[\mr{g}_{\ast2}\big]
+\mathbf{Te}^{2}\,\big[-\frac{5}{2}\,\zeta(5)\big]\,\mr{g}_{\ast2}^4
+\mathbf{Te}^{4}\,\big[\frac{1}{4}\,\zeta(3)\big]\,\mr{g}_{\ast2}^4
+\mathbf{Te}^{3}\,\big[\frac{3}{4}\,\zeta(3)^2\big]\,\mr{g}_{\ast2}^5
\\
+\mathbf{Te}^{2}
  \!\!\!&\!\!&\!\!\!
\big[
\frac{3}{2}\,\zeta(3)^3
-\frac{4}{5}\,\zeta(3)\,\zeta(2)^3
+\frac{47}{6}\,\zeta(9)\big]\,\mr{g}_{\ast2}^6
+\mathbf{Te}^{4}\,\big[
-\frac{21}{40}\,\zeta(3)\,\zeta(2)^2
-\frac{63}{64}\,\zeta(7)\big]\,\mr{g}_{\ast2}^6
\\
+\mathbf{Te}^{6}
  \!\!\!&\!\!&\!\!\!
\big[
\frac{3}{8}\,\zeta(2)\,\zeta(3)
+\frac{1}{16}\,\zeta(5)\big]\,\mr{g}_{\ast2}^6
+\mathbf{Te}^{8}\,\big[
-\frac{1}{16}\,\zeta(3)\big]\,\mr{g}_{\ast2}^6
\end{eqnarray*}

\noindent 
{\bf Example 5}:\, $\mg{\gop_\ast}$ up to weight 15 for $f=l\circ g$ with 
$ g(z)=z\,\Big[1+3\,g_{\ast3} z^{-3} \Big]^{\frac{1}{3}}$.
\begin{eqnarray*}
+\mathbf{Te}^{2}
   \!\!\!&\!\!&\!\!\!
\big[ \mr{g}_{\ast3}\big]
+\mathbf{Te}^{2} \big[ \m 6\,\zeta(3)\big]\,\mr{g}_{\ast3}^{2}
+\mathbf{Te}^{2} 
\big[\m \frac{32}{5}\,\zeta(2)^3+36\,\zeta(3)^2\big]\,\mr{g}_{\ast3}^{3}
+\mathbf{Te}^{4} \big[\m \frac{1}{5}\,\zeta(2)^2\big]\,\mr{g}_{\ast3}^{3}
\\
+\mathbf{Te}^{6} 
   \!\!\!&\!\!&\!\!\!
\big[\m \frac{1}{3}\,\zeta(2)\big]\,\mr{g}_{\ast3}^{3}
+\mathbf{Te}^{2}
 \big[\m 210\,\zeta(9)\p\frac{576}{5}\,\zeta(3)\,\zeta(2)^3\m 216\,\zeta(3)^3\big]\,\mr{g}_{\ast3}^{4}
 \\
+\mathbf{Te}^{4} 
   \!\!\!&\!\!&\!\!\!
\big[14\,\zeta(7)\p \frac{18}{5}\,\zeta(3)\,\zeta(2)^2\big]\,\mr{g}_{\ast3}^{4}
+\mathbf{Te}^{6} \big[6\,\zeta(2)\,\zeta(3)\m \frac{10}{3}\,\zeta(5)\big]\,\mr{g}_{\ast3}^{4}
\\
+\mathbf{Te}^{2} 
   \!\!\!&\!\!&\!\!\!
\big[1960\,\zeta(7)\,\zeta(5)\p 5880\,\zeta(3)\,\zeta(9)\m \frac{6912}{5}\,\zeta(3)^2\,\zeta(2)^3   
          \m \frac{23054144}{125125}\,\zeta(2)^6
\\
   \!\!\!&\!\!&\!\!\!
             \p 1296\,\zeta(3)^4\m 420\,\zeta(10,2)\big]\,\mr{g}_{\ast3}^{5}
             \\
+\mathbf{Te}^{4}
                \!\!\!&\!\!&\!\!\!
 \big[\m\frac{216}{5}\,\zeta(3)^2\,\zeta(2)^2\m 434\,\zeta(3)\,\zeta{7}]
+\frac{1332224}{28875}\,\zeta(2)^5\m 38\,\zeta(5)^2\m 49\,\zeta(8,2)\big]\,\mr{g}_{\ast3}^{5}
\\
+\mathbf{Te}^{6} 
   \!\!\!&\!\!&\!\!\!
\big[
\m 72\,\zeta(2)\,\zeta(3)^2\p\frac{340}{3}\,\zeta(3)\,\zeta(5)
\p\frac{1007}{1575}\,\zeta(2)^4
\m\frac{50}{3}\,\zeta(6,2)\big]\,\mr{g}_{\ast3}^{5}
\\ 
+\mathbf{Te}^{8}
   \!\!\!&\!\!&\!\!\!
 \big[
\frac{193}{75}\,\zeta(2)^3\big]\,\mr{g}_{\ast3}^{5}
+\mathbf{Te}^{10} \big[
\frac{16}{15}\,\zeta(2)^2\big]\,\mr{g}_{\ast3}^{5}
+\mathbf{Te}^{12} \big[
\frac{7}{45}\,\zeta(2) \big]\,\mr{g}_{\ast3}^{5}
\end{eqnarray*}





%

{\bf REFERENCES.}
\\
\\ \noindent
{\bf [B]} O. Bouillot, {\it Invariants analytiques des diff\'{e}omorphismes et multiz\^{e}tas},
Orsay PhD, 19.10.2011, available as PDF on O. Bouillot's homepage.
\\ \noindent
{\bf [BEE]} X. Buff, J. Ecalle, A. Epstein, {\it Limits of degenerate parabolic quadratic rational maps},
\\ \noindent
{\bf [E1]} J.Ecalle, {\it Alg\`ebres de fonctions r\'esurgentes.} Pub. Math. Orsay (1981).
\\
{\bf [E2]} J.Ecalle,  {\it Les fonctions r\'esurgentes appliqu\'ees \`a l'it\' eration.} Pub. Math. Orsay (1981).
\\ \noindent
{\bf [E3]} J. Ecalle,
{\it L'\'equation du pont et la classification analytique des objets locaux.}
Pub. Math. Orsay (1985).
\\ \noindent
{\bf [E4]} J.  Ecalle, {\it   Twisted Resurgence Monomials and canonical synthesis of Local
Objects.} Proc. of the June 2002 Edinburgh conference
 on Analysable Functions (ICMS workshop), World Scientific.
\\  \noindent
{\bf [F] } P. Fatou, {\it Sur les solutions uniformes de certaines \'equations fonctionnelles.} C. R. Acad. Sci. Paris, 
{\bf 143}, 546-548 (1906).




%
\end{document}